# IDENTIFICATION OF MULTITYPE BRANCHING PROCESSES

By F. Maaouia and A. Touati

*Faculté des Sciences de Tunis and Faculté des Sciences de Bizerte*

We solve the problem of constructing an asymptotic global confidence region for the means and the covariance matrices of the reproduction distributions involved in a supercritical multitype branching process. Our approach is based on a central limit theorem associated with a quadratic law of large numbers performed by the maximum likelihood or the multidimensional Lotka–Nagaev estimator of the reproduction law means. The extension of this approach to the least squares estimator of the mean matrix is also briefly discussed.

On résout le problème de construction d'une région de confiance asymptotique et globale pour les moyennes et les matrices de covariance des lois de reproduction d'un processus de branchement multitype et supercritique. Notre approche est basée sur un théorème de limite centrale associé à une loi forte quadratique vérifiée par l'estimateur du maximum de vraisemblance ou l'estimateur multidimensionnel de Lotka–Nagaev des moyennes des lois de reproduction. L'extension de cette approche à l'estimateur des moindres carrés de la matrice des moyennes est aussi brièvement commentée.

## 1. Introduction.

1.1. *Motivation.* Statistical inference about the means and/or the covariance matrices of the reproduction distributions involved in a positively regular supercritical *Bienaymé–Galton–Watson* process with $d$-types [BGW($d$)] has been investigated by several authors [2, 4, 5, 14, 26]. Though some important work (discussed below) has been done on this topic, a satisfactory global approach has not been outlined. The purpose of this article is to fill this gap.

For the convenience of the reader, our method is initially derived in the familiar context of the one-type BGW process. So, let $\mathbf{X} = (\mathbf{X}_n)_{n \geq 0}$ be a supercritical process starting from $\mathbf{X}_0 = 1$. The identification of the mean $a = \mathbb{E}(\mathbf{X}_1) \in \,]1, +\infty[$ and the variance $\sigma^2 = \mathrm{Var}(\mathbf{X}_1) \in \,]0, +\infty[$ of the offspring distribution of $\mathbf{X}$ is a classical problem which has been studied by many. An exhaustive review of this topic can be found in [14].

Let us recall some significant results in connection with this subject. The maximum likelihood estimator of $a$, given by

$$(1.1) \qquad \widehat{a}_n = (\mathbf{S}_{n-1})^{-1}(\mathbf{S}_n - 1), \qquad \mathbf{S}_n = \sum_{k=0}^{n} \mathbf{X}_k,$$

satisfies the relation

$$(1.2) \quad \mathbf{S}_{n-1}(\widehat{a}_n - a) = \sum_{k=1}^{n}(\mathbf{X}_k - a\mathbf{X}_{k-1}) = \sum_{k=1}^{n}(\mathbf{X}_k - \mathbb{E}(\mathbf{X}_k/\mathcal{F}_{k-1})) = \mathbf{L}_n,$$

where $\mathcal{F}_n$ is the $\sigma$-algebra generated by the random variables (r.v.'s) $\mathbf{X}_0, \ldots, \mathbf{X}_n$. It is strongly consistent on the set of nonextinction $\mathbf{E} = \{\lim_{n \to \infty} \mathbf{S}_n = \infty\} = \{\lim_{n \to \infty} \mathbf{X}_n = \infty\}$. Moreover, conditional on $\mathbf{E}$, the central limit theorem (CLT) with random normalization,

$$(1.3) \qquad\qquad \sqrt{\mathbf{S}_{n-1}}(\widehat{a}_n - a) \xrightarrow[n \to \infty]{\mathcal{L}} \mathcal{N}(0, \sigma^2),$$

holds [11, 12, 13, 14, 16, 17, 19, 20, 21, 23, 24]; $\mathcal{N}(0, \sigma^2)$ denotes the centered Gaussian distribution with variance $\sigma^2$. As pointed out by Dion [13, 14], conditional on the set $\mathbf{E}_n = \{\mathbf{X}_n > 0\}$, the result (1.3) is also true.

To build a confidence region for $a$, it remains to estimate $\sigma^2$. Dion [13] and Heyde [18] proved that

$$
\begin{aligned}
(1.4) \quad \check{\sigma}_n^2 &= \frac{1}{n}\sum_{k=1}^{n}\mathbf{X}_{k-1}^{-1}(\mathbf{X}_k - \check{a}_n\mathbf{X}_{k-1})^2 \\
&= \frac{1}{n}\sum_{k=1}^{n}\mathbf{X}_{k-1}(\check{a}_k - \check{a}_n)^2 \xrightarrow[n \to \infty]{} \sigma^2 \qquad \text{a.s. on } \mathbf{E},
\end{aligned}
$$

where $\check{a}_n = \mathbf{X}_{n-1}^{-1}\mathbf{X}_n \mathbf{1}_{\{\mathbf{X}_{n-1} > 0\}}$ is the *Lotka–Nagaev estimator* of $a$, subsequently called the *empirical estimator* of $a$, which is also strongly consistent on $\mathbf{E}$. In addition, under the assumption $\mathbb{E}(\mathbf{X}_1^4) < \infty$, they showed that

$$\sqrt{n}(\check{\sigma}_n^2 - \sigma^2) \xrightarrow[n \to \infty]{\mathcal{L}} \mathcal{N}(0, 2\sigma^4).$$

We improve their result here by showing that, conditional on the set $\mathbf{E}$ or $\mathbf{E}_n$,

$$(1.5) \quad \{\sqrt{n}(\sigma^{-2}\check{\sigma}_n^2 - 1), \sqrt{\mathbf{S}_{n-1}}\check{\sigma}_n^{-1}(\widehat{a}_n - a)\} \xrightarrow[n \to \infty]{\mathcal{L}} \mathcal{N}(0, 2) \otimes \mathcal{N}(0, 1)$$

and

$$(1.6) \qquad \limsup \sqrt{\frac{n}{\ln \ln n}} |\sigma^{-2} \check{\sigma}_n^2 - 1| = 2 \qquad \text{a.s.,}$$

where $\otimes$ denote the tensor product of measures.

A similar approach for estimating $\sigma^2$ may be derived from the following property satisfied by the martingale $(\mathbf{L}_n)$, called the *quadratic strong law of large numbers* (QSL) in [8, 9] (see also [28]):

$$(1.7) \qquad \frac{1}{n} \sum_{k=1}^{n} (\mathbf{S}_{k-1}^{-1/2} \mathbf{L}_k)^2 \underset{n \to \infty}{\longrightarrow} \sigma^2 \qquad \text{a.s. on } \mathbf{E}.$$

From (1.2) and (1.7), it is easy to derive that

$$(1.8) \qquad \widehat{\sigma}_n^2 = \frac{1}{n} \sum_{k=1}^{n} \mathbf{S}_{k-1} (\widehat{a}_k - \widehat{a}_n)^2 \underset{n \to \infty}{\longrightarrow} \sigma^2 \qquad \text{a.s. on } \mathbf{E}.$$

We shall prove here that, under the additional hypothesis $\mathbb{E}(\mathbf{X}_1^4) < \infty$ and conditional on the set $\mathbf{E}$, the following two properties hold:

$$(1.9) \qquad \{\sqrt{n}(\sigma^{-2}\widehat{\sigma}_n^2 - 1), \sqrt{\mathbf{S}_{n-1}} \widehat{\sigma}_n^{-1} (\widehat{a}_n - a)\} \underset{n \to \infty}{\overset{\mathcal{L}}{\longrightarrow}} \mathcal{N}\left(0, 2\frac{a+1}{a-1}\right) \otimes \mathcal{N}(0, 1);$$

$$(1.10) \qquad \limsup \sqrt{\frac{n}{\ln \ln n}} |\sigma^{-2} \widehat{\sigma}_n^2 - 1| = 2\sqrt{\frac{a+1}{a-1}} \qquad \text{a.s.}$$

Hence, the estimator $\check{\sigma}_n^2$ is asymptotically more efficient than $\widehat{\sigma}_n^2$; however, it is insensitive to any change that occurs on the mean.

In the remainder of this article, the global approach we develop for the one-type BGW process is generalized to the $d$-type case. Moreover, results analogous to (1.8), (1.9) and (1.10) for the empirical estimator and the least squares estimator of the reproduction law means are discussed.

## 1.2. *Assumptions.*

### 1.2.1. *About the observed sample.*

From now on $^*\mathbf{X}_n = (\mathbf{X}_n(1), \ldots, \mathbf{X}_n(d))$ denotes the generic state of a BGW($d$) process, that is, the column vector of numbers of particles (or individuals) of each type in the $n$th generation. The initial state $\mathbf{X}_0$ is taken equal to the vector $\mathbf{1}$ whose components are equal to 1.

The particles of type $j$ that are alive in the $(n-1)$st generation give birth to a total number of $\mathbf{Y}_n^j(1)$ particles of type $1, \ldots, \mathbf{Y}_n^j(d)$ particles of type $d$. Therefore, we have

$$(1.11) \qquad \mathbf{X}_n = \sum_{j=1}^{d} \mathbf{Y}_n^j, \qquad \mathbf{Y}_n^j = {}^*(\mathbf{Y}_n^j(1), \ldots, \mathbf{Y}_n^j(d)).$$

Except in Section 5, we suppose that the r.v.'s $(\mathbf{X}_0, \mathbf{Y}_1^j, \ldots, \mathbf{Y}_n^j; 1 \leq j \leq d)$ are observable and we denote by $\mathcal{G}_n$ the $\sigma$-algebra they generate.

1.2.2. *About the reproduction laws.* The following common assumptions will be used subsequently.

ASSUMPTION A-1. For each $j \in \{1, \ldots, d\}$, the reproduction distribution $\nu_j$ of the particle (or individual) of type $j$ is assumed to have finite moments of order $2\delta$ ($\delta \in \,]1, 2]$) and an invertible covariance matrix.

ASSUMPTION A-2. The BGW($d$) process $(\mathbf{X}_n)_{n \geq 0}$ is assumed to be positively regular and supercritical (see [3, 25] for the details).

The column vector that represents the mean of $\nu_j$ and its covariance matrix are denoted $a^j$ and $K^j$, respectively.

1.3. *Brief description of the main results.* In the three next sections we are concerned with estimation of the parameters $\{a^j, K^j\}_{1 \leq j \leq d}$ based on observation of the sample $\{\mathbf{X}_0, \mathbf{Y}_k^j; 1 \leq k \leq n, 1 \leq j \leq d\}$. More precisely, we analyze in Section 2 the asymptotic behavior of the *maximum likelihood estimator* (MLE) $\widehat{\mathcal{A}}_n$ and the *empirical estimator* (EE) $\widecheck{\mathcal{A}}_n$ of $\mathcal{A} = \text{Vect}(a^1, \ldots, a^d)$ (see Section 1.4 below). For each one of these estimators, we shall (1) prove strong consistency on the set of nonextinction $\mathbf{E}$ and also give the strong rate of convergence, and (2) prove asymptotic normality conditional on the set $\mathbf{E}_n = \bigcup_{j=1}^d \{\mathbf{X}_n(j) > 0\}$ after appropriate centering and random normalization.

In Section 3 we show that $\widehat{\mathcal{A}}_n$ and $\widecheck{\mathcal{A}}_n$ satisfy two QSLs as in the case of the one-type BGW process. This allows us to derive two estimators $\widehat{\mathcal{K}}_n$ and $\widecheck{\mathcal{K}}_n$ of $\mathcal{K} = \text{Diag}\{K^1, \ldots, K^d\}$. We prove also their strong consistency on $\mathbf{E}$ and sharpen their strong rate of convergence.

The construction of a global confidence region for the parameters $\{a^j, K^j\}_{1 \leq j \leq d}$ is achieved in Section 4 thanks to a central limit theorem performed by the pair $(\widehat{\mathcal{A}}_n, \widehat{\mathcal{K}}_n)$ or $(\widecheck{\mathcal{A}}_n, \widecheck{\mathcal{K}}_n)$.

Finally, in Section 5, a similar approach is discussed for the *least squares estimator* (LSE) of the matrix $\mathbf{A} = [a^1, \ldots, a^d]$.

1.4. *Notation.* We recall some standard notation.

N1. $\mathbf{I}_d$ and $\mathcal{I}_{d^2}$ denote, respectively, the $d \times d$ and the $d^2 \times d^2$ identity matrices.

N2. For a real $d \times d$ matrix $B$ whose column vectors are $b^1, \ldots, b^d$,

$$\text{Vect}(B) = \text{Vect}(b^1, \ldots, b^d)$$
$$= {}^*(b^1(1), \ldots, b^1(d), \ldots, b^d(1), \ldots, b^d(d)) \in \mathbb{R}^{d^2}.$$

N3. The notation $\mathrm{Diag}(B_1, B_2)$ refers to the block matrix $\begin{bmatrix} B_1 & 0 \\ 0 & B_2 \end{bmatrix}$.

N4. If $A = (a_{i,j})_{i,j}$ and $B = (b_{r,l})_{r,l}$ are two matrices, $A \otimes B = (a_{i,j}B)_{i,j}$ designates the block matrix whose $(i, j)$ block is $a_{i,j}B$.

We also introduce the following less common notation on block matrices:

N5. If $P = (P(i, j))_{1 \le i,j \le d}$ and $Q = (Q(r, l))_{1 \le r,l \le d}$ are two block matrices, then $P \boxtimes Q = ((P(i, j) \otimes Q(r, l))_{1 \le r,l \le d})_{1 \le i,j \le d}$ stands for the block matrix whose $(i, j)$ block is the block matrix $(P(i, j) \otimes Q(r, l))_{1 \le r,l \le d}$.

N6. If $C = (C(i, j))_{1 \le i,j \le d}$ is a block matrix, then $^{\perp}C$ denotes the block matrix $(C(j, i))_{1 \le i,j \le d}$.

The standard Euclidean inner product on $\mathbb{R}^q$ and its associated norm are denoted by $\langle \cdot, \cdot \rangle$ and $\| \cdot \|$, respectively. The trace operator on square matrices is denoted by $\mathrm{tr}(\cdot)$.

## 2. Estimation of the reproduction means.
From the definition of $\mathbf{Y}_n^j$,

$$(2.1) \quad \widehat{a}_n^j = \left( \sum_{p=0}^{n-1} \mathbf{X}_p(j) \right)^{-1} \sum_{p=1}^{n} \mathbf{Y}_p^j \quad \text{and} \quad \check{a}_n^j = \mathbf{X}_{n-1}(j)^{-1} \mathbf{Y}_n^j \mathbf{1}_{\{\mathbf{X}_{n-1}(j) > 0\}}$$

are intuitively candidates for estimating $a^j$. Clearly, $\check{a}_n^j$ is the multidimensional analog of the empirical estimator defined in Section 1.1. As is shown below, $\widehat{a}_n^j$ is in fact the MLE of $a^j$, that is, the estimator derived by maximizing each component of $a^j$.

2.1. *Computation of the maximum likelihood estimator of $\mathcal{A}$.* Let $\mathcal{F}_n$ denote the $\sigma$-algebra generated by the r.v.'s $(\mathbf{X}_0, \mathbf{X}_1, \dots, \mathbf{X}_n)$. The conditional distribution of $\mathbf{X}_{n+1}$ with respect to $\mathcal{F}_n$ is $\nu_1^{*\mathbf{X}_n(1)} * \cdots * \nu_d^{*\mathbf{X}_n(d)}$, where $\nu_j^{*k}$ represents the convolution product of $k$ distributions equal to $\nu_j$. Consequently, the following useful construction of the Markov chain $(\mathbf{X}_n)$ is available. Given a probability space $(\Omega, \mathcal{F}, \mathbb{P})$, let $\{(\xi_{n,k}^j)_{(n,k) \in \mathbb{N}^2}\}_{1 \le j \le d}$ be a system of $d$ independent sequences of r.v.'s on $(\Omega, \mathcal{F}, \mathbb{P})$, where $(\xi_{n,k}^j)_{(n,k) \in \mathbb{N}^2}$ are i.i.d. random vectors with respect to the distribution $\nu_j$. More precisely, the vector $\xi_{n,k}^j$ represents the successors of the $k$th particle (or individual) of type $j$ that are alive in the $(n-1)$st generation. Hence, we have

$$(2.2) \quad \mathbf{Y}_n^j = \left( \sum_{k=1}^{\mathbf{X}_{n-1}(j)} \xi_{n,k}^j \right) \mathbf{1}_{\{\mathbf{X}_{n-1}(j) \ne 0\}}, \qquad j \in \{1, \dots, d\}.$$

The likelihood function of the r.v.'s $\{\mathbf{X}_0, \dots, \mathbf{X}_n\}$ is

$$\mathfrak{L}_n = \mathfrak{L}_n(\mathbf{X}_0, \dots, \mathbf{X}_n) = \prod_{p=1}^{n} \nu_1^{*\mathbf{X}_{p-1}(1)} * \cdots * \nu_d^{*\mathbf{X}_{p-1}(d)}(\mathbf{X}_p)$$

(2.3)
$$= \prod_{p=1}^{n} \left( \prod_{j=1}^{d} \left( \prod_{l \in \mathcal{D}_j} \nu_j(l)^{\mathbf{N}_p^j(l)} \right) \right),$$

where for each $j \in \{1, \ldots, d\}$ the r.v.'s $\mathbf{N}_p^j(l) = \sum_{k=1}^{\mathbf{X}_{p-1}(j)} \mathbf{1}_{\{\xi_{p,k}^j = l\}}, \ l \in \mathbb{N}^d$, satisfy the properties

$$\sum_{l \in \mathbb{N}^d} \mathbf{N}_p^j(l) = \mathbf{X}_{p-1}(j), \qquad \sum_{l \in \mathbb{N}^d} l \mathbf{N}_p^j(l) = \mathbf{Y}_p^j,$$

and for $1 \leq j \leq d$, $\mathcal{D}_j = \{l \in \mathbb{N}^d / \nu_j(l) > 0\}$ is the support of the reproduction distribution $\nu_j$. Therefore, using the Lagrangian technique, we prove that the r.v.'s

(2.4)
$$\widehat{\nu_j(l)}_n = \left( \sum_{p=0}^{n-1} \mathbf{X}_p(j) \right)^{-1} \sum_{p=1}^{n} \mathbf{N}_p^j(l), \qquad l \in \mathcal{D}_j,$$

maximize $\mathfrak{L}_n$, but we emphasize the fact that $\widehat{\nu_j(l)}_n$ is an estimator of $\nu_j(l)$ only if the r.v.'s $(\xi_{n,k}^j)$ are observable. However, the MLE of the mean $a^j$ of $\nu_j$, obtained by maximizing each component of $a^j$, that is,

(2.5)
$$\sum_{l \in \mathbb{N}^d} l \widehat{\nu_j(l)}_n = \left( \sum_{p=0}^{n-1} \mathbf{X}_p(j) \right)^{-1} \sum_{p=1}^{n} \mathbf{Y}_p^j,$$

coincides with $\widehat{a}_n^j$. It is a statistic of the observed sample $\{\mathbf{X}_0, \mathbf{Y}_1^j, \ldots, \mathbf{Y}_n^j, 1 \leq j \leq d\}$. Hence, the MLE of $\mathcal{A}$ is

(2.6)
$$\widehat{\mathcal{A}}_n = \text{Vect}(\widehat{a}_n^1, \ldots, \widehat{a}_n^d) = \mathcal{S}_{n-1}^{-1} \sum_{p=1}^{n} \text{Vect}(\mathbf{Y}_p^1, \ldots, \mathbf{Y}_p^d),$$

where $\mathcal{S}_n = \text{Diag}(\mathbf{S}_n(1)\mathbf{I}_d, \ldots, \mathbf{S}_n(d)\mathbf{I}_d)$ and $\mathbf{S}_{n-1}(j) = \sum_{p=0}^{n-1} \mathbf{X}_p(j)$.

In particular, let us remark that

(2.7)
$$\mathcal{S}_{n-1}(\widehat{\mathcal{A}}_n - \mathcal{A}) = \sum_{k=1}^{n} \zeta_k, \qquad \text{where } \zeta_k = \text{Vect}(\zeta_k^1, \ldots, \zeta_k^d),$$

and

$$\zeta_k^j = \mathbf{Y}_k^j - \mathbb{E}(\mathbf{Y}_k^j / \mathcal{G}_{k-1}) = \left( \sum_{r=1}^{\mathbf{X}_{k-1}(j)} (\xi_{k,r}^j - a^j) \right) \mathbf{1}_{\{\mathbf{X}_{k-1}(j) > 0\}}.$$

2.2. *Asymptotic properties of the MLE of $\mathcal{A}$.*

2.2.1. *Main result.* The asymptotic properties of $(\widehat{\mathcal{A}}_n)$, announced in Section 1.3 and stated below, are in fact direct consequences of those of the martingale defined by the right-hand side of (2.7).

THEOREM 2.1. *Let* $\mathbf{X} = (\Omega, \mathcal{F}, \mathbb{P}, (\mathbf{X}_n)_{n \geq 0})$ *be a d-multitype branching process, starting from* $\mathbf{1} = {}^*(1, \ldots, 1)$, *whose reproduction distributions satisfy Assumptions* A-1 *and* A-2. *Then the following properties hold, where* $\mathbf{E}$ *designates the set of nonextinction of* $\mathbf{X}$ *and* $\mathbf{E}_n = \bigcup_{j=1}^{d} \{\mathbf{X}_n(j) > 0\}$.

(i) *The estimator* $(\widehat{\mathcal{A}}_n)$ *is a strongly consistent estimator of* $\mathcal{A}$ *on the set* $\mathbf{E}$; *more precisely,*

$$\max_{1 \leq k \leq n} \|\rho^{k/2}(\widehat{\mathcal{A}}_k - \mathcal{A})\|^2 = O(\ln n) \qquad a.s. \ on \ \mathbf{E}.$$

(ii) *The estimator* $(\widehat{\mathcal{A}}_n)$ *is asymptotically normal; more precisely, conditional on the set* $\mathbf{E}$ *or* $\mathbf{E}_n$,

$$\mathcal{S}_{n-1}^{1/2}(\widehat{\mathcal{A}}_n - \mathcal{A}) \xrightarrow[n \to \infty]{\mathcal{L}} \mathcal{N}_{d^2}(0, \mathcal{K}).$$

REMARK 2.2. The estimation of the offspring means of a supercritical multitype branching process has been studied by Asmussen and Keiding [2], Keiding and Lauritzen [22] and Nanthi [1, 26]. Their results are strictly contained in Theorem 2.1.

The proof of Theorem 2.1 is based on some important asymptotic properties of the martingale defined by (2.7), which is stated in Lemma 2.3 below.

2.2.2. *Auxiliary results for the proof of Theorem* 2.1. Before stating the lemma, let us recall that Assumption A-2 allows us to affirm that the following three conditions are satisfied by $\mathbf{A} = [a^1, \ldots, a^d]$ (the matrix whose column vectors are $a^1, \ldots, a^d$):

C1. The matrix $\mathbf{A}$ has a maximal eigenvalue $\rho > 1$, which is equal to the spectral radius of $\mathbf{A}$. Moreover, the modulus of each other eigenvalue of $\mathbf{A}$ is strictly less than $\rho$.

C2. There exist an eigenvector $u = {}^*(u(1), \ldots, u(d))$ of $\mathbf{A}$ and an eigenvector $v = {}^*(v(1), \ldots, v(d))$ of ${}^*\mathbf{A}$ (the transpose matrix of $\mathbf{A}$), associated to $\rho$, such that

$$u(j) > 0, \ \forall 1 \leq j \leq d, \qquad \langle v, \mathbf{1} \rangle = \sum_{j=1}^{d} v(j) = 1,$$

$$\langle u, v \rangle = \sum_{j=1}^{d} u(j)v(j) = 1.$$

C3. If $\mathbf{P} = u^*v$, then there exits a matrix $\mathbf{R}$ that satisfies

$$\mathbf{PR} = \mathbf{RP}, \qquad \mathbf{R}^n = O(n^{d-1}\rho_0^n) \qquad \text{for } \rho_0 \in ]0, \rho[, \qquad \mathbf{A}^n = \rho^n \mathbf{P} + \mathbf{R}^n.$$

Conditions C1 and C2 imply that the set of nonextinction of the BGW($d$) process $(\mathbf{X}_n)$, that is, $\mathbf{E} = \limsup \mathbf{E}_n$, where $\mathbf{E}_n = \bigcup_{j=1}^d \{\mathbf{X}_n(j) > 0\}$, is not negligible. In fact, one can prove the existence of a r.v. $\mathbf{W}$ such that $\mathbf{E} = \{\mathbf{W} > 0\}$, $\mathbb{P}(\mathbf{E}) > 0$ and the following properties hold almost surely:

(P$_1$)  $\rho^{-n}\langle v, \mathbf{X}_n \rangle \underset{n \to \infty}{\longrightarrow} \mathbf{W}$, $\rho^{-n}\mathbf{X}_n \underset{n \to \infty}{\longrightarrow} \mathbf{W}u$.

Moreover, $\mathbb{P}(\mathbf{W} = 0) = 0$ if for all $1 \le j \le d$ the reproduction distribution $\nu_j$ belongs to a regular exponential model.

Using the obvious recursive relation

$$\mathbf{X}_n = \mathbf{A}\mathbf{X}_{n-1} + \varepsilon_n, \qquad \varepsilon_n = \sum_{j=1}^d \zeta_n^j,$$

condition C3, combined with the law of the iterated logarithm, allows us to prove the following property more precisely than (P$_1$):

(P$_1'$)  $\rho^{-n}\mathbf{X}_n - \mathbf{W}u = O(\sqrt{\ln n}\, n^{d-1}\theta^n)$ a.s. on $\mathbf{E}$, where $\theta = (\frac{\max(\rho_0, 1)}{\rho})$.

Let $(\mathcal{M}_n)$ be the normalized martingale defined from the right-hand side of (2.7) by

$$
\begin{aligned}
\mathcal{M}_n &= \mathcal{U}^{-1/2}\mathcal{K}^{-1/2}\sum_{k=1}^n \zeta_k = \text{Vect}(\mathbf{M}_n^1, \ldots, \mathbf{M}_n^d), \\
(2.8) & \\
\mathbf{M}_n^j &= u(j)^{-1/2}(K^j)^{-1/2}\sum_{k=1}^n \zeta_k^j, \qquad \mathcal{U} = \text{Diag}(u(1)\mathbf{I}_d, \ldots, u(d)\mathbf{I}_d).
\end{aligned}
$$

Note that its predictable quadratic variation is

$$
\begin{aligned}
(2.9) \quad \langle \mathcal{M} \rangle_n &= \text{Diag}(u(1)^{-1}\mathbf{S}_{n-1}(1)\mathbf{I}_d, \ldots, u(d)^{-1}\mathbf{S}_{n-1}(d)\mathbf{I}_d) \\
&= \mathcal{U}^{-1/2}\mathcal{S}_{n-1}\mathcal{U}^{-1/2}.
\end{aligned}
$$

LEMMA 2.3. *The following properties hold for the martingale $(\mathcal{M}_n)$, where $\mathbf{E}$, $\mathbf{E}_n$ are the sets defined above in Theorem* 2.1.

(P$_2$)  *We have $\rho^{-n}\langle \mathcal{M} \rangle_n \underset{n \to \infty}{\longrightarrow} \frac{\mathbf{W}}{\rho - 1}\mathcal{I}_{d^2}$ a.s. on $\mathbf{E}$, where $\mathbf{W}$ is the r.v. defined in* (P$_1$).

(P$_3$)  *We also have $\max_{1 \le k \le n} \|\rho^{-k/2}\mathcal{M}_k\|^2 = O(\ln n)$ a.s. on $\mathbf{E}$.*

(P$_4$) *For* $t \in \mathbb{R}^{d^2}$ *and* $k \leq n$, *let* $\varphi_{n,k}(t) = \mathbb{E}(\exp\{i\langle t, \rho^{-n/2}\Delta\mathcal{M}_k\rangle\}|\mathcal{G}_{k-1})$ *and* $\Phi_n(t) = \prod_{k=1}^{n}\varphi_{n,k}(t)$, *where* $\Delta\mathcal{M}_k = \mathcal{M}_k - \mathcal{M}_{k-1}$. *Then*

$$\Phi_n(t) \underset{n\to\infty}{\longrightarrow} \exp\left\{-\frac{1}{2}\frac{\mathbf{W}}{\rho-1}\|t\|^2\right\} \qquad a.s.\ on\ \mathbf{E}.$$

(P$_5$) *Conditional on* $\mathbf{E}$ *or* $\mathbf{E}_n$, $\rho^{-n/2}\mathcal{M}_n \overset{\mathcal{L}}{\longrightarrow}_{n\to\infty} \Sigma((\mathbf{W}/\rho-1)^{1/2}\mathcal{I}_{d^2})$ *stably, where* $\Sigma(\mathcal{T})$ *is a r.v. independent of* $\mathbf{W}$ *with the same distribution as* $\mathcal{N}_{d^2}(0, \mathcal{T}^*\mathcal{T})$.

PROOF OF THEOREM 2.1.   To get the strong consistency of the estimator $(\widehat{\mathcal{A}}_n)$ on $\mathbf{E}$, we note from (2.7), (2.8) and (2.9) that

$$\rho^{n/2}(\widehat{\mathcal{A}}_n - \mathcal{A}) = \rho^{n/2}\mathcal{S}_{n-1}^{-1}\mathcal{U}^{1/2}\mathcal{K}^{1/2}\mathcal{M}_n$$
$$= \mathcal{U}^{-1/2}(\rho^{-n}\langle\mathcal{M}\rangle_n)^{-1}\mathcal{K}^{1/2}(\rho^{-n/2}\mathcal{M}_n),$$

but thanks to properties (P$_2$) and (P$_3$) satisfied by the martingale $\mathcal{M}$, we have

$$\max_{1\leq k\leq n}\|\rho^{-k/2}\mathcal{M}_k\|^2 = O(\ln n) \quad \text{and} \quad \rho^{-n}\langle\mathcal{M}\rangle_n \underset{n\to\infty}{\longrightarrow} \frac{\mathbf{W}}{\rho-1}\mathcal{I}_{d^2} \qquad \text{a.s. on } \mathbf{E}.$$

Hence, (i) is proved.

By (P$_1$), we have $\rho^{-n}\mathcal{S}_{n-1} \longrightarrow_{n\to\infty} \frac{\mathbf{W}}{\rho-1}\mathcal{U}$ a.s., so property (ii) of the theorem is a consequence of (P$_5$), since

$$(2.10) \qquad \mathcal{S}_{n-1}^{1/2}(\widehat{\mathcal{A}}_n - \mathcal{A}) = (\rho^{n/2}\mathcal{S}_{n-1}^{-1/2})\mathcal{U}^{1/2}\mathcal{K}^{1/2}(\rho^{-n/2}\mathcal{M}_n). \qquad \square$$

PROOF OF LEMMA 2.3.   Equality (2.9) and property (P$_1$) imply (P$_2$), that is,

$$\rho^{-n}\langle\mathcal{M}\rangle_n = \mathcal{U}^{-1/2}(\rho^{-n}\mathcal{S}_{n-1})\mathcal{U}^{-1/2} \underset{n\to\infty}{\longrightarrow} \frac{\mathbf{W}}{\rho-1}\mathcal{I}_{d^2} \qquad \text{a.s. on } \mathbf{E}.$$

Property (P$_3$) is obtained thanks to the following result satisfied by the r.v.'s $\zeta_n = \mathrm{Vect}(\zeta_n^1, \ldots, \zeta_n^d)$:

$$(2.11) \qquad \begin{aligned} \zeta_n^j &= \sum_{k=1}^{\mathbf{X}_{n-1}(j)}(\xi_{n,k}^j - a^j) \\ &= O(\sqrt{\mathbf{X}_{n-1}(j)\ln\ln\mathbf{X}_{n-1}(j)}) \\ &= O(\sqrt{\rho^n \ln n}) \qquad \text{on } \mathbf{E}. \end{aligned}$$

In fact, (2.11) is a consequence of (P$_1$) and the LIL property

$$\sum_{k=1}^{\mathbf{N}}(\xi_{n,k}^j - a^j) = O(\sqrt{2\mathbf{N}\ln\ln\mathbf{N}}) \qquad \text{a.s.}$$

Since we have

$$\mathcal{M}_n = \mathcal{U}^{-1/2}\mathcal{K}^{-1/2}\left(\sum_{k=1}^{n}\zeta_k\right)$$

$$= O\left(\left(\sum_{k=1}^{n}\rho^{k/2}\right)\sqrt{\ln n}\right)$$

$$= O(\rho^{n/2}\sqrt{\ln n}),$$

then

$$\max_{1\leq k\leq n}\|\rho^{-k/2}\mathcal{M}_k\|^2 = O(\ln n) \qquad \text{a.s. on } \mathbf{E}.$$

Let us now prove properties (P$_4$) and (P$_5$). For $t_1,\ldots,t_d\in\mathbb{R}^d$, $k\leq n$ and $t=\text{Vect}(t_1,\ldots,t_d)$, we have

$$\varphi_{n,k}(t) = \mathbb{E}\left(\prod_{j=1}^{d}\exp\left\{i\left\langle t_j, \frac{1}{\sqrt{\rho^n u(j)}}(K^j)^{-1/2}\zeta_k^j\right\rangle\right\}\Big|\mathcal{G}_{k-1}\right),$$

so

$$\Phi_n(t) = \prod_{k=1}^{n}\varphi_{n,k}(t)$$

$$= \prod_{j=1}^{d}\left[\mathbb{E}\left(\exp\left\{i\left\langle t_j, \frac{1}{\sqrt{\rho^n u(j)}}(K^j)^{-1/2}(\xi_{1,1}^j - a^j)\right\rangle\right\}\right)\right]^{\mathbf{S}_{n-1}(j)}$$

$$\times \mathbf{1}_{\bigcap_{k=1}^{n}\{\mathbf{X}_{k-1}(j)>0\}}.$$

By noting that

$$\Phi_n(t) \sim \prod_{j=1}^{r}\exp\left\{-\frac{1}{2}\frac{1}{\rho^n u(j)}\mathbf{S}_{n-1}(j)\|t_j\|^2\right\}\mathbf{1}_{\bigcap_{k=1}^{n}\{\mathbf{X}_{k-1}(j)>0\}} \qquad \text{a.s. on } \mathbf{E},$$

we deduce that

$$\Phi_n(t) \underset{n\to\infty}{\longrightarrow} \prod_{j=1}^{d}\exp\left\{-\frac{1}{2}\frac{\mathbf{W}}{\rho-1}\|t_j\|^2\right\} = \exp\left\{-\frac{1}{2}\frac{\mathbf{W}}{\rho-1}\|t\|^2\right\} \qquad \text{a.s. on } \mathbf{E}.$$

Property (P$_2$) combined with property (P$_4$) allows us to use Theorem 3 of [27], which implies that, conditional on the set $\mathbf{E} = \{\mathbf{W}>0\}$,

$$(\rho^{-n/2}\mathcal{M}_n) \underset{n\to\infty}{\overset{\mathcal{L}}{\longrightarrow}} \Sigma\left(\left(\frac{\mathbf{W}}{\rho-1}\right)^{1/2}\mathcal{I}_{d^2}\right) \qquad \text{stably,}$$

where $\Sigma(\mathcal{T})$ is a r.v. independent of $\mathbf{W}$ with the same distribution as $\mathcal{N}_{d^2}(0,\mathcal{T}^*\mathcal{T})$. As in [12], we can affirm that the above result is also true conditional on $\mathbf{E}_n$. $\quad\square$

2.3. *Asymptotic properties of the empirical estimator of $\mathcal{A}$.* From the definition (2.1) of $\check{a}_n^j$, we deduce that for $j = 1, 2, \ldots, d$,

$$\mathbf{X}_{n-1}(j)(\check{a}_n^j - a^j) = \sum_{k=1}^{\mathbf{X}_{n-1}(j)} (\xi_{n,k}^j - a^j) \mathbf{1}_{\{\mathbf{X}_{n-1}(j) > 0\}} = \zeta_n^j.$$

Hence, setting $\mathcal{X}_n = \text{Diag}(\mathbf{X}_n(1)\mathbf{I}_d, \ldots, \mathbf{X}_n(d)\mathbf{I}_d)$, the empirical estimator $\check{\mathcal{A}}_n$ of $\mathcal{A}$ satisfies the relations

$$(2.12) \quad \mathcal{X}_{n-1}(\check{\mathcal{A}}_n - \mathcal{A}) = \zeta_n \quad \text{and} \quad \mathcal{S}_{n-1}(\widehat{\mathcal{A}}_n - \mathcal{A}) = \sum_{k=1}^{n} \mathcal{X}_{k-1}(\check{\mathcal{A}}_k - \mathcal{A}).$$

So property (i) of Theorem 2.1 holds for $\check{\mathcal{A}}_n$. An analogous result of property (ii) of this theorem may also be stated.

## 3. Estimation of the reproduction covariance matrices.

3.1. *Some asymptotic properties of the sequences $(\mathcal{S}_{n-1}^{1/2}(\widehat{\mathcal{A}}_n - \mathcal{A}))$ and $(\mathcal{X}_{n-1}^{1/2}(\check{\mathcal{A}}_n - \mathcal{A}))$.* As announced in Section 1.3, the key tool for constructing a strong estimator of the covariance matrices $\mathcal{K} = \text{Diag}(K^1, \ldots, K^d)$ is the QSL that underlies the CLT satisfied by the sequence $(\mathcal{S}_{n-1}^{1/2}(\widehat{\mathcal{A}}_n - \mathcal{A}))$ or $(\mathcal{X}_{n-1}^{1/2}(\check{\mathcal{A}}_n - \mathcal{A}))$. This property is stated in the second part of the next theorem.

THEOREM 3.1. *Within the framework of Theorem 2.1, the sequence $(\mathcal{S}_{n-1}^{1/2}(\widehat{\mathcal{A}}_n - \mathcal{A}))$ or $(\mathcal{X}_{n-1}^{1/2}(\check{\mathcal{A}}_n - \mathcal{A}))$ satisfies the following properties almost surely on $\mathbf{E}$.*

(iii) *ASCLT:*

$$\frac{1}{n} \sum_{k=1}^{n} \delta_{\mathcal{S}_{k-1}^{1/2}(\widehat{\mathcal{A}}_k - \mathcal{A})} \underset{n \to \infty}{\Longrightarrow} \mathcal{N}_{d^2}(0, \mathcal{K}),$$

$$\frac{1}{n} \sum_{k=1}^{n} \delta_{\mathcal{X}_{k-1}^{1/2}(\check{\mathcal{A}}_k - \mathcal{A})} \underset{n \to \infty}{\Longrightarrow} \mathcal{N}_{d^2}(0, \mathcal{K})$$

($\Longrightarrow$ *denotes the weak convergence of measures*).

(iv) *QSL:*

$$\frac{1}{n} \sum_{k=1}^{n} \mathcal{S}_{k-1}^{1/2}(\widehat{\mathcal{A}}_k - \mathcal{A})^*(\widehat{\mathcal{A}}_k - \mathcal{A})\mathcal{S}_{k-1}^{1/2} \underset{n \to \infty}{\longrightarrow} \mathcal{K},$$

$$\frac{1}{n} \sum_{k=1}^{n} \mathcal{X}_{k-1}^{1/2}(\check{\mathcal{A}}_k - \mathcal{A})^*(\check{\mathcal{A}}_k - \mathcal{A})\mathcal{X}_{k-1}^{1/2} \underset{n \to \infty}{\longrightarrow} \mathcal{K}.$$

*Moreover, if the reproduction distributions have finite moments of order* 4, *then*

(v) *LIL of the QSL:*

$$\limsup \sqrt{\left(\frac{n}{\ln \ln n}\right)} \left| \frac{1}{n} \sum_{k=1}^{n} \|\mathcal{K}^{-1/2} \mathcal{S}_{k-1}^{1/2}(\widehat{\mathcal{A}}_k - \mathcal{A})\|^2 - d^2 \right| = 2d\sqrt{\frac{\rho+1}{\rho-1}}$$

*and*

$$\limsup \sqrt{\left(\frac{n}{\ln \ln n}\right)} \left| \frac{1}{n} \sum_{k=1}^{n} \|\mathcal{K}^{-1/2} \mathcal{X}_{k-1}^{1/2}(\widecheck{\mathcal{A}}_k - \mathcal{A})\|^2 - d^2 \right| = 2d.$$

3.2. *Construction of two strong estimators of* $\mathcal{K}$. Using property (iv), we construct two strong consistent estimators of $\mathcal{K}$ as stated in the theorem below.

**THEOREM 3.2.** *Within the framework of Theorem* 2.1, *let*

$$\widehat{K}_n^j = \frac{1}{n} \sum_{k=1}^{n} \mathrm{S}_{k-1}(j)(\widehat{a}_k^j - \widehat{a}_n^j)^*(\widehat{a}_k^j - \widehat{a}_n^j), \qquad j = 1, \dots, d \quad and$$

$$\widehat{\mathcal{K}}_n = \mathrm{Diag}(\widehat{K}_n^1, \dots, \widehat{K}_n^d);$$

$$\widecheck{K}_n^j = \frac{1}{n} \sum_{k=1}^{n} \mathrm{X}_{k-1}(j)(\widecheck{a}_k^j - \widecheck{a}_n^j)^*(\widecheck{a}_k^j - \widecheck{a}_n^j), \qquad j = 1, \dots, d \quad and$$

$$\widecheck{\mathcal{K}}_n = \mathrm{Diag}(\widecheck{K}_n^1, \dots, \widecheck{K}_n^d).$$

(vi) *Then on the set of nonextinction* **E**,

$$\widehat{\mathcal{K}}_n \underset{n \to \infty}{\longrightarrow} \mathcal{K} \qquad a.s., \qquad \widecheck{\mathcal{K}}_n \underset{n \to \infty}{\longrightarrow} \mathcal{K} \qquad a.s.$$

(vii) *Moreover, if the reproduction distributions have finite moments of order* 4, *then*

$$\limsup \sqrt{\frac{n}{\ln \ln n}} |\mathrm{tr}(\mathcal{K}^{-1/2}\widehat{\mathcal{K}}_n \mathcal{K}^{-1/2}) - d^2| = 2d\sqrt{\frac{\rho+1}{\rho-1}} \qquad a.s.,$$

$$\limsup \sqrt{\frac{n}{\ln \ln n}} |\mathrm{tr}(\mathcal{K}^{-1/2}\widecheck{\mathcal{K}}_n \mathcal{K}^{-1/2}) - d^2| = 2d \qquad a.s.$$

**REMARK 3.3.** In [26], Nanthi proposed an estimator of a generic element of $K^j$ by adapting the Dion [13] and Heyde [18] method for estimating the variance of the offspring distribution of a BGW process (see Section 1.1).

Let us first prove Theorem 3.2.

PROOF OF THEOREM 3.2. We only prove the results announced for the sequence $(\mathcal{S}_{n-1}^{1/2}(\widehat{\mathcal{A}}_n - \mathcal{A}))$, because those corresponding to the sequence $(\mathcal{X}_{n-1}^{1/2}(\widecheck{\mathcal{A}}_n - \mathcal{A}))$ can be established in the same way.

By the property (iv) of Theorem 3.1, we can affirm that, for all $1 \leq j \leq d$,

$$\frac{1}{n}\sum_{k=1}^{n}\mathbf{S}_{k-1}(j)(\widehat{a}_k^j - a^j)^*(\widehat{a}_k^j - a^j) \underset{n\to\infty}{\longrightarrow} K^j \qquad \text{a.s. on } \mathbf{E}.$$

Since almost surely on $\mathbf{E}$

$$\max_{1\leq k\leq n}\|\mathbf{S}_{k-1}(j)^{1/2}(\widehat{a}_k^j - a^j)\|^2 = O\left(\max_{1\leq k\leq n}\|\mathcal{S}_{k-1}^{1/2}(\widehat{\mathcal{A}}_k - \mathcal{A})\|^2\right) = O(\ln n)$$

and

$$\sum_{k=1}^{n}\mathbf{S}_{k-1}(j) = O(\mathbf{S}_{n-1}(j))$$

both hold, we deduce that almost surely on $\mathbf{E}$,

$$\frac{1}{n}\left(\sum_{k=1}^{n}S_{k-1}(j)\right)\|\widehat{a}_n^j - a^j\|^2 = O\left(\frac{\ln n}{n}\right)$$

and

$$\frac{1}{n}\left(\sum_{k=1}^{n}S_{k-1}(j)\|\widehat{a}_k^j - a^j\|\right)\|\widehat{a}_n^j - a^j\|$$

$$\leq \frac{1}{n}\left(\sum_{k=1}^{n}\sqrt{S_{k-1}(j)}\max_{1\leq r\leq n}\sqrt{S_{r-1}(j)}\|\widehat{a}_r^j - a^j\|\right)\|\widehat{a}_n^j - a^j\| = O\left(\frac{\ln n}{n}\right).$$

However, the estimator $\widehat{K}_n^j$ satisfies the equality

$$\widehat{K}_n^j - \frac{1}{n}\sum_{k=1}^{n}\mathbf{S}_{k-1}(j)(\widehat{a}_k^j - a^j)^*(\widehat{a}_k^j - a^j)$$

$$= -\frac{1}{n}\left(\sum_{k=1}^{n}\mathbf{S}_{k-1}(j)(\widehat{a}_k^j - a^j)\right)^*(\widehat{a}_n^j - a^j)$$

$$-\frac{1}{n}(\widehat{a}_n^j - a^j)\left(\sum_{k=1}^{n}\mathbf{S}_{k-1}(j)^*(\widehat{a}_k^j - a^j)\right)$$

$$+\frac{1}{n}\left(\sum_{k=1}^{n}\mathbf{S}_{k-1}(j)\right)(\widehat{a}_n^j - a^j)^*(\widehat{a}_n^j - a^j);$$

(viii) hence

$$\widehat{\mathcal{K}}_n - \mathcal{K} = \frac{1}{n} \sum_{k=1}^{n} \mathcal{S}_{k-1}^{1/2} (\widehat{\mathcal{A}}_k - \mathcal{A})^* (\widehat{\mathcal{A}}_k - \mathcal{A}) \mathcal{S}_{k-1}^{1/2} - \mathcal{K} + O\left(\frac{\ln n}{n}\right) \qquad \text{a.s. on } \mathbf{E}.$$

This property shows the strong consistency of $\widehat{\mathcal{K}}_n$ on $\mathbf{E}$. Moreover, combined with property (v) of Theorem 3.1, it allows us to affirm that assertion (vii) of Theorem 3.2 is also true for $\widehat{\mathcal{K}}_n$.   $\square$

The proof of Theorem 3.1 becomes simple thanks to the next results, which establish that properties (iii), (iv) and (v) are in fact the transcriptions of analogous ones satisfied by the normalized martingale $(\rho^{-n/2}\mathcal{M}_n)$ or the sequence $(\zeta_n)$ related to $(\mathcal{M}_n)$ by

(3.1)          $$\zeta_n = \text{Vect}(\zeta_n^1, \ldots, \zeta_n^d), \qquad \zeta_n^j = \sqrt{u(j)}(K^j)^{1/2} \Delta \mathbf{M}_n^j.$$

Note that

(3.2)          $$\mathbb{E}(\zeta_n/\mathcal{G}_{n-1}) = 0 \quad \text{and} \quad \mathbb{E}(\zeta_n^*\zeta_n/\mathcal{G}_{n-1}) = \mathcal{X}_{n-1}^{1/2}\mathcal{K}\mathcal{X}_{n-1}^{1/2}.$$

### 3.3. Auxiliary results for the proof of Theorem 3.1.

LEMMA 3.4.  *For the sequences* $(\mathcal{M}_n)$ *and* $(\zeta_n)$, *defined, respectively, by* (2.8) *and* (3.1), *the following properties hold almost surely on* $\mathbf{E}$.

$(\mathrm{P}_6)$  *ASCLT:*

$$\frac{1}{n} \sum_{k=1}^{n} \delta_{\rho^{-k/2}\mathcal{M}_k} \quad \underset{n\to\infty}{\Longrightarrow} \quad \Sigma\left(\left(\frac{\mathbf{W}}{\rho-1}\right)^{1/2}\mathcal{I}_{d^2}\right),$$

$$\frac{1}{n} \sum_{k=1}^{n} \delta_{\mathcal{X}_{k-1}^{-1/2}\zeta_k} \quad \underset{n\to\infty}{\Longrightarrow} \quad \Sigma(\mathcal{K}),$$

*where* $\Sigma(\mathcal{X})$ *is a Gaussian r.v. as in property (*$\mathrm{P}_5$*).*

$(\mathrm{P}_7)$  *QSL:*

$$\frac{1}{n} \sum_{k=1}^{n} \rho^{-k}(\mathcal{M}_k^*\mathcal{M}_k - \langle\mathcal{M}\rangle_k) \underset{n\to\infty}{\longrightarrow} 0,$$

$$\text{so } \frac{1}{n} \sum_{k=1}^{n} \rho^{-k}\mathcal{M}_k^*\mathcal{M}_k \underset{n\to\infty}{\longrightarrow} \frac{\mathbf{W}}{\rho-1}\mathcal{I}_{d^2},$$

$$\frac{1}{n} \sum_{k=1}^{n} \mathcal{X}_{k-1}^{-1/2}(\zeta_k^*\zeta_k - \mathcal{X}_{k-1}^{1/2}\mathcal{K}\mathcal{X}_{k-1}^{1/2})\mathcal{X}_{k-1}^{-1/2} \underset{n\to\infty}{\longrightarrow} 0,$$

$$\text{so } \frac{1}{n} \sum_{k=1}^{n} \mathcal{X}_{k-1}^{-1/2}\zeta_k^*\zeta_k\mathcal{X}_{k-1}^{-1/2} \underset{n\to\infty}{\longrightarrow} \mathcal{K}.$$

*Moreover, if the reproduction distributions have finite moments of order 4, then:*

(P$_8$) *LIL of the QSL:*

$$\limsup \sqrt{\left(\frac{n}{\ln\ln n}\right)} \left|\frac{1}{n}\sum_{k=1}^{n} \rho^{-k}\operatorname{tr}(\mathcal{M}_k{}^*\mathcal{M}_k - \langle\mathcal{M}\rangle_k)\right| = \frac{(\rho+1)^{1/2}2\,d\mathbf{W}}{(\rho-1)^{3/2}},$$

$$\limsup \sqrt{\left(\frac{n}{\ln\ln n}\right)} \left|\frac{1}{n}\sum_{k=1}^{n} \operatorname{tr}(\mathcal{K}^{-1/2}\mathcal{X}_{k-1}^{-1/2}(\zeta_k{}^*\zeta_k - \mathcal{X}_{k-1}^{1/2}\mathcal{K}\mathcal{X}_{k-1}^{1/2})\right.$$
$$\left. \times\,\mathcal{X}_{k-1}^{-1/2}\mathcal{K}^{-1/2})\right| = 2d.$$

The proof of this lemma, except property (P$_8$), is postponed to Appendix A. Property (P$_8$) is proved at the end of Appendix B.

LEMMA 3.5. *Within the framework of the Theorem* 2.1:

(ix) *The properties*

$$\sum_{k=1}^{n}\mathcal{S}_{k-1}^{1/2}(\widehat{\mathcal{A}}_k - \mathcal{A})^*(\widehat{\mathcal{A}}_k - \mathcal{A})\mathcal{S}_{k-1}^{1/2}$$
$$= n\mathcal{K} + \left(\frac{\rho-1}{\mathbf{W}}\right)\mathcal{K}^{1/2}\left(\sum_{k=1}^{n}\rho^{-k}(\mathcal{M}_k{}^*\mathcal{M}_k - \langle\mathcal{M}\rangle_k)\right)\mathcal{K}^{1/2} + o(n^{1/2}),$$

$$\sum_{k=1}^{n}\mathcal{X}_{k-1}^{1/2}(\widecheck{\mathcal{A}}_k - \mathcal{A})^*(\widecheck{\mathcal{A}}_k - \mathcal{A})\mathcal{X}_{k-1}^{1/2}$$
$$= n\mathcal{K} + \sum_{k=1}^{n}\mathcal{X}_{k-1}^{-1/2}(\zeta_k{}^*\zeta_k - \mathcal{X}_{k-1}^{1/2}\mathcal{K}\mathcal{X}_{k-1}^{1/2})\mathcal{X}_{k-1}^{-1/2}$$

*hold almost surely on* $\mathbf{E}$.

(x) *Consequently,*

$$\sum_{k=1}^{n}\|\mathcal{K}^{-1/2}\mathcal{S}_{k-1}^{1/2}(\widehat{\mathcal{A}}_k - \mathcal{A})\|^2 - nd^2$$
$$= \left(\frac{\rho-1}{\mathbf{W}}\right)\sum_{k=1}^{n}\rho^{-k}\operatorname{tr}(\mathcal{M}_k{}^*\mathcal{M}_k - \langle\mathcal{M}\rangle_k) + o(n^{1/2}),$$

$$\sum_{k=1}^{n}\|\mathcal{K}^{-1/2}\mathcal{X}_{k-1}^{1/2}(\widecheck{\mathcal{A}}_k - \mathcal{A})\|^2 - nd^2$$
$$= \sum_{k=1}^{n}\operatorname{tr}(\mathcal{K}^{-1/2}\mathcal{X}_{k-1}^{-1/2}(\zeta_k{}^*\zeta_k)\mathcal{X}_{k-1}^{-1/2}\mathcal{K}^{-1/2} - \mathcal{I}_{d^2}).$$

This lemma is proved in Appendix A.

Proof of Theorem 3.1.    We only prove the results announced for the martingale $(\mathcal{M}_n)$ because the results that correspond to the sequence $(\zeta_n)$ can be obtained in the same way.

The property $\rho^{-n}\langle\mathcal{M}\rangle_n \longrightarrow_{n\to\infty} \frac{\mathbf{W}}{\rho-1}\,\mathcal{I}_{d^2}$ a.s. on $\mathbf{E}$, combined with the ASCLT satisfied by $\mathcal{M}$ [see (P$_6$)], implies that, conditional on the set of nonextinction $\mathbf{E}$,

$$\frac{1}{n}\sum_{k=1}^{n}\delta_{\left(\rho^{-k/2}\mathcal{M}_k,\rho^{-k}\langle\mathcal{M}\rangle_k\right)} \underset{n\to\infty}{\Longrightarrow} \mu\otimes\delta_{\mathcal{C}} \qquad \text{a.s.,}$$

where $\mu$ is the law of the r.v. $\Sigma\left(\left(\frac{\mathbf{W}}{\rho-1}\right)^{1/2}\mathcal{I}_{d^2}\right)$ and $\mathcal{C}=\left(\frac{\mathbf{W}}{\rho-1}\right)\mathcal{I}_{d^2}$ is the almost-sure limit of $(\rho^{-k}\langle\mathcal{M}\rangle_k)$ on $\mathbf{E}$. Hence we deduce property (iii) of the theorem for $(\mathcal{M}_n)$ thanks to (2.7), (2.8) and (2.9).

Property (iv) for $(\mathcal{M}_n)$ is a consequence of the QSL satisfied by the martingale $\mathcal{M}$ [see (P$_7$)] and the relation (ix) of Lemma 3.5.

Likewise, the LIL stated in the theorem, for the sequence $(\mathcal{S}_{n-1}(\widehat{\mathcal{A}}_n-\mathcal{A}))$, is a consequence of the LIL satisfied by $(\mathcal{M}_n)$ [see (P$_8$)] and the relation (x) of Lemma 3.5.   □

## 4. Asymptotic confidence region for the parameters $\{a^j, K^j\}_j$.

Our goal here is to construct an asymptotic confidence region for all the parameters $\{a^j, K^j\}_{1\le j\le d}$. The key tool is the CLT stated below for the pair of estimators $(\widehat{\mathcal{A}}_n, \widehat{\mathcal{K}}_n)$ and $(\widehat{\mathcal{A}}_n, \overleftarrow{\mathbb{K}}_n)$.

4.1. *Central limit theorem for the pair of estimators* $(\widehat{\mathcal{A}}_n, \widehat{\mathcal{K}}_n)$ *and* $(\widehat{\mathcal{A}}_n, \overleftarrow{\mathbb{K}}_n)$. The proof of the next theorem, based on some probabilistic results performed by the sequences $(\mathcal{M}_n)$ and $(\zeta_n)$, is postponed to the end of Section 4.

Theorem 4.1.    *Let* $\mathbf{X}=(\Omega,\mathcal{F},\mathbb{P},\mathbb{F},(\mathbf{X}_n)_{n\ge 0})$ *be a $d$-multitype branching process, starting from* $\mathbf{1}={}^*(1,\ldots,1)$. *We assume that its reproduction distributions satisfy Assumption A-1 with* $\delta=2$ *and also Assumption A-2. Then, conditional on the set* $\mathbf{E}$ *or* $\mathbf{E}_n$:

(xi)  *We have*

$$\{\sqrt{n}(\widehat{\mathcal{K}}_n-\mathcal{K}),\mathcal{S}_{n-1}^{1/2}(\widehat{\mathcal{A}}_n-\mathcal{A})\} \xrightarrow[n\to\infty]{\mathcal{L}} \left\{\mathfrak{G}_1\left(\sqrt{\frac{2}{\rho-1}}\mathcal{K}\right)+\mathfrak{G}_2(\mathcal{K}),\mathfrak{G}\right\},$$

$$\{\sqrt{n}(\overleftarrow{\mathbb{K}}_n-\mathcal{K}),\mathcal{S}_{n-1}^{1/2}(\widehat{\mathcal{A}}_n-\mathcal{A})\} \xrightarrow[n\to\infty]{\mathcal{L}} \{\mathfrak{G}_2(\mathcal{K}),\mathfrak{G}\},$$

*where* $\mathfrak{G}$ *is a r.v. distributed as* $\mathcal{N}_{d^2}(0,\mathcal{K})$ *and, for* $\mathcal{T}=\mathrm{Diag}(T^1,\ldots,T^d)$, $(\mathfrak{G}_r(\mathcal{T}))_{r\in\{1,2\}}$ *are independent identically distributed Gaussian matrices,* *which are also independent of* $\mathfrak{G}$. *Moreover the covariance matrix of* $\mathrm{Vect}(\mathfrak{G}_r(\mathcal{T}))$ *is equal to*

$$\mathrm{Diag}(T^j \otimes T^j + ^{\perp}(\mathrm{Vect}(T^j)^* \mathrm{Vect}(T^j)); 1 \leq j \leq d).$$

(xii) *We also have*

$$\{\sqrt{n}(\mathrm{tr}(\mathcal{K}^{-1/2}\widehat{\mathcal{K}}_n\mathcal{K}^{-1/2}) - d^2), \mathcal{S}_{n-1}^{1/2}(\widehat{\mathcal{A}}_n - \mathcal{A})\}$$

$$\xrightarrow[n\to\infty]{\mathcal{L}} \mathcal{N}\left(0, 2d^2\frac{\rho+1}{\rho-1}\right) \otimes \mathcal{N}_{d^2}(0, \mathcal{K}),$$

$$\{\sqrt{n}(\mathrm{tr}(\mathcal{K}^{-1/2}\check{\mathcal{K}}_n\mathcal{K}^{-1/2}) - d^2), \mathcal{S}_{n-1}^{1/2}(\widehat{\mathcal{A}}_n - \mathcal{A})\}$$

$$\xrightarrow[n\to\infty]{\mathcal{L}} \mathcal{N}(0, 2d^2) \otimes \mathcal{N}_{d^2}(0, \mathcal{K}).$$

4.2. *Construction of confidence regions.* From property (xii), and the fact that

$$(4.1) \qquad \frac{\langle \mathbf{X}_n, \mathbf{1}\rangle}{\langle \mathbf{X}_n, \mathbf{1}\rangle + 2\langle \mathbf{S}_{n-1}, \mathbf{1}\rangle} \xrightarrow[n\to\infty]{} \frac{\rho-1}{\rho+1} \qquad \text{a.s. on } \mathbf{E},$$

we deduce the following result, which allows us to construct an asymptotic confidence region for all the parameters $\{a^j, K^j\}_{1 \leq j \leq d}$.

THEOREM 4.2. *Conditional on the set* $\mathbf{E}$ *or* $\mathbf{E}_n$:

(xiii) *We have*

$$\left\{\sqrt{\frac{n\langle \mathbf{X}_n, \mathbf{1}\rangle}{2d^2(\langle \mathbf{X}_n, \mathbf{1}\rangle + 2\langle \mathbf{S}_{n-1}, \mathbf{1}\rangle)}}\left(\sum_{j=1}^d \mathrm{tr}((K^j)^{-1/2}\widehat{K}_n^j(K^j)^{-1/2}) - d^2\right),\right.$$

$$\left.\sum_{j=1}^d \mathbf{S}_{k-1}(j)\|(\widehat{K}_n^j)^{-1/2}(\widehat{a}_k^j - a^j)\|^2\right\}$$

$$\xrightarrow[n\to\infty]{\mathcal{L}} \mathcal{N}(0, 1) \otimes \chi(d^2),$$

$$\left\{\sqrt{\frac{n}{2d^2}}\left(\sum_{j=1}^d \mathrm{tr}((K^j)^{-1/2}\check{K}_n^j(K^j)^{-1/2}) - d^2\right),\right.$$

$$\left.\sum_{j=1}^d \mathbf{S}_{k-1}(j)\|(\check{K}_n^j)^{-1/2}(\widehat{a}_k^j - a^j)\|^2\right\}$$

$$\xrightarrow[n\to\infty]{\mathcal{L}} \mathcal{N}(0, 1) \otimes \chi(d^2),$$

*where* $\chi(d^2)$ *denote the chi-square distribution with* $d^2$ *degrees of freedom.*

REMARK 4.3. Let $\{a_0^j, K_0^j\}_{1 \le j \le d}$ be a given possible structure of the reproduction means and covariance matrices. From Theorem 4.2 one can easily derive an asymptotic test of the hypothesis $\{a^j, K^j\}_{1 \le j \le d} = \{a_0^j, K_0^j\}_{1 \le j \le d}$ against the alternative $\{a^j, K^j\}_{1 \le j \le d} \ne \{a_0^j, K_0^j\}_{1 \le j \le d}$.

## 4.3. Auxiliary results for the proof of Theorem 4.1.

4.3.1. *The CLT of the QSL satisfied by $(\mathcal{M}_n)$ or $(\zeta_n)$.* For the QSL satisfied by the martingale $(\mathcal{M}_n)$ or the sequence $(\zeta_n)$, that is, property $(P_7)$, the rate of the weak convergence is given by the following technical lemma, which will be established in Appendix B.

LEMMA 4.4. *If the assumptions of Theorem 4.1 are satisfied, then conditional on the set $\mathbf{E}$ or $\mathbf{E}_n$:*

$(P_9)$ *We have*

$$\left\{ \frac{1}{\sqrt{n}} \sum_{k=1}^n \rho^{-k} (\mathcal{M}_k{}^* \mathcal{M}_k - \langle \mathcal{M} \rangle_k), \rho^{-n/2} \mathcal{M}_n \right\}$$

$$\xrightarrow[n \to \infty]{\mathcal{L}} \left\{ \Sigma_1 \left( \left( \frac{2\mathbf{W}^2}{(\rho-1)^3} \right)^{1/4} \mathbf{I}_d \right) + \Sigma_2 \left( \sqrt{\frac{\mathbf{W}}{\rho-1}} \mathbf{I}_d \right), \Sigma \left( \left( \frac{\mathbf{W}}{\rho-1} \right)^{1/2} \mathcal{I}_{d^2} \right) \right\},$$

$$\left\{ \frac{1}{\sqrt{n}} \sum_{k=1}^n \mathcal{K}^{-1/2} \mathcal{X}_{k-1}^{-1/2} (\zeta_k{}^* \zeta_k - \mathcal{X}_{k-1}^{1/2} \mathcal{K} \mathcal{X}_{k-1}^{1/2}) \mathcal{X}_{k-1}^{-1/2} \mathcal{K}^{-1/2}, \rho^{-n/2} \mathcal{M}_n \right\}$$

$$\xrightarrow[n \to \infty]{\mathcal{L}} \left\{ \Sigma_2(\mathbf{I}_d), \Sigma \left( \left( \frac{\mathbf{W}}{\rho-1} \right)^{1/2} \mathcal{I}_{d^2} \right) \right\},$$

*where $\mathbf{W}, \Sigma(\cdot)$ are r.v.'s as in property $(P_5)$ and $(\Sigma_r(T))_{r \in \{1,2\}}$ are independent identically distributed Gaussian matrices, which are also independent of the pair $(\mathbf{W}, \Sigma(\cdot))$. The covariance matrix of $\mathrm{Vect}(\Sigma_r(T))$ is equal to*

$$(T \otimes T) \otimes (T \otimes T) + {}^{\perp}(\mathrm{Vect}\, T^* \mathrm{Vect}\, T) \boxtimes {}^{\perp}(\mathrm{Vect}\, T^* \mathrm{Vect}\, T).$$

$(P'_9)$ *In particular,*

$$\left\{ \left[ \frac{1}{\sqrt{n}} \sum_{k=1}^n \rho^{-k} (\mathbf{M}_k^{j}{}^* \mathbf{M}_k^j - \langle \mathbf{M}^j \rangle_k) \right]_{1 \le j \le d}, \rho^{-n/2} \mathcal{M}_n \right\}$$

$$\xrightarrow[n \to \infty]{\mathcal{L}} \left\{ \left[ \Sigma_1^j \left( \frac{\sqrt{2}\mathbf{W}}{(\rho-1)^{3/2}} \mathbf{I}_d \right) + \Sigma_2^j \left( \frac{\mathbf{W}}{\rho-1} \mathbf{I}_d \right) \right]_{1 \le j \le d},$$

$$\Sigma \left( \left( \frac{\mathbf{W}}{\rho-1} \right)^{1/2} \mathcal{I}_{d^2} \right) \right\}$$

$$\left\{\left[\frac{1}{\sqrt{n}}\sum_{k=1}^{n}\mathbf{X}_{k-1}(j)^{-1}(K^j)^{-1/2}(\zeta_k^{j*}\zeta_k^j - X_{k-1}(j)K^j)(K^j)^{-1/2}\right]_{1\le j\le d},\right.$$
$$\left.\rho^{-n/2}\mathcal{M}_n\right\}$$

$$\xrightarrow[n\to\infty]{\mathcal{L}}\left\{[\Sigma_2^j(\mathbf{I}_d)]_{1\le j\le d}, \Sigma\left(\left(\frac{\mathbf{W}}{\rho-1}\right)^{1/2}\mathcal{I}_{d^2}\right)\right\},$$

where $\{\Sigma_r^j(T), 1\le j\le d, r\in\{1,2\}\}$ are independent identically distributed Gaussian matrices, which are also independent of the pair $(\mathbf{W},\Sigma(\cdot))$. The covariance matrix of $\mathrm{Vect}(\Sigma_r^j(T))$ is equal to $T\otimes T + {}^{\perp}(\mathrm{Vect}\,T^*\mathrm{Vect}\,T)$.

(P$_{10}$) Moreover,

$$\left\{\frac{1}{\sqrt{n}}\sum_{k=1}^{n}\rho^{-k}\,\mathrm{tr}(\mathcal{M}_k{}^*\mathcal{M}_k - \langle\mathcal{M}\rangle_k), \rho^{-n/2}\mathcal{M}_n\right\}$$

$$\xrightarrow[n\to\infty]{\mathcal{L}}\left\{\left(\frac{2\mathbf{W}^2 d^2(\rho+1)}{(\rho-1)^3}\right)^{1/2}G, \Sigma\left(\left(\frac{\mathbf{W}}{\rho-1}\right)^{1/2}\mathcal{I}_{d^2}\right)\right\},$$

$$\left\{\frac{1}{\sqrt{n}}\sum_{k=1}^{n}\mathrm{tr}(\mathcal{K}^{-1/2}\mathcal{X}_{k-1}^{-1/2}(\zeta_k{}^*\zeta_k - \mathcal{X}_{k-1}^{1/2}\mathcal{K}\mathcal{X}_{k-1}^{1/2})\mathcal{X}_{k-1}^{-1/2}\mathcal{K}^{-1/2}), \rho^{-n/2}\mathcal{M}_n\right\}$$

$$\xrightarrow[n\to\infty]{\mathcal{L}}\left\{\sqrt{2}dG, \Sigma\left(\left(\frac{\mathbf{W}}{\rho-1}\right)^{1/2}\mathcal{I}_{d^2}\right)\right\},$$

where $G$ is a standard Gaussian r.v., which is independent of the pair $\{\mathbf{W},\Sigma(\cdot)\}$.

4.3.2. *The CLT of the QSL satisfied by* $(\widehat{\mathcal{A}}_n)$ *and* $(\widecheck{\mathcal{A}}_n)$. Combined with Lemma 3.5, the previous result allows us to establish the following one.

LEMMA 4.5. *If the assumptions of Theorem 4.1 are satisfied, then conditional on* $\mathbf{E}$ *or* $\mathbf{E}_n$:

(xiv) *We have*

$$\left\{\sqrt{n(\rho-1)}\left(\frac{1}{n}\sum_{k=1}^{n}\mathcal{K}^{-1/2}\mathcal{S}_{k-1}^{1/2}(\widehat{\mathcal{A}}_k - \mathcal{A})^*(\widehat{\mathcal{A}}_k - \mathcal{A})\mathcal{S}_{k-1}^{1/2}\mathcal{K}^{-1/2} - \mathcal{I}_{d^2}\right),\right.$$
$$\left.\mathcal{S}_{n-1}^{1/2}(\widehat{\mathcal{A}}_n - \mathcal{A})\right\}$$

$$\xrightarrow[n\to\infty]{\mathcal{L}}\{\Sigma_1(2^{1/4}\mathbf{I}_d) + \Sigma_2((\rho-1)^{1/4}\mathbf{I}_d), \mathfrak{G}\},$$

$$\left\{\sqrt{n}\left(\frac{1}{n}\sum_{k=1}^{n}\mathcal{K}^{-1/2}\mathcal{X}_{k-1}^{1/2}(\widecheck{\mathcal{A}}_k - \mathcal{A})^*(\widecheck{\mathcal{A}}_k - \mathcal{A})\mathcal{X}_{k-1}^{1/2}\mathcal{K}^{-1/2} - \mathcal{I}_{d^2}\right),\right.$$

$$\mathcal{S}_{n-1}^{1/2}(\widehat{\mathcal{A}}_n - \mathcal{A})\Big\}$$

$$\xrightarrow[n\to\infty]{\mathcal{L}} \{\Sigma_2(\mathbf{I}_d), \mathfrak{G}\},$$

*where $\mathfrak{G}$ is a r.v. distributed as $\mathcal{N}_{d^2}(0, \mathcal{K})$ and $(\Sigma_r(T))_{r\in\{1,2\}}$ are independent identically distributed Gaussian matrices, which are also independent of $\mathfrak{G}$. The covariance matrix of $\mathrm{Vect}(\Sigma_r(T))$ is equal to*

$$(T \otimes T) \otimes (T \otimes T) + {}^\perp(\mathrm{Vect}\, T^* \mathrm{Vect}\, T) \boxtimes {}^\perp(\mathrm{Vect}\, T^* \mathrm{Vect}\, T).$$

(xv) *We also have*

$$\left\{ \sqrt{\frac{n(\rho-1)}{(\rho+1)}} \left( \frac{1}{n} \sum_{k=1}^n \|\mathcal{K}^{-1/2} \mathcal{S}_{k-1}^{1/2}(\widehat{\mathcal{A}}_k - \mathcal{A})\|^2 - d^2 \right), \mathcal{S}_{n-1}^{1/2}(\widehat{\mathcal{A}}_n - \mathcal{A}) \right\}$$

$$\xrightarrow[n\to\infty]{\mathcal{L}} \mathcal{N}(0, 2d^2) \otimes \mathcal{N}_{d^2}(0, \mathcal{K}),$$

$$\left\{ \sqrt{n} \left( \frac{1}{n} \sum_{k=1}^n \|\mathcal{K}^{-1/2} \mathcal{X}_{k-1}^{1/2}(\widehat{\mathcal{A}}_k - \mathcal{A})\|^2 - d^2 \right), \mathcal{S}_{n-1}^{1/2}(\widehat{\mathcal{A}}_n - \mathcal{A}) \right\}$$

$$\xrightarrow[n\to\infty]{\mathcal{L}} \mathcal{N}(0, 2d^2) \otimes \mathcal{N}_{d^2}(0, \mathcal{K}).$$

PROOF OF THEOREM 4.1.    To prove property (xi) of the theorem for the sequence $(\widehat{\mathcal{K}}_n)$, first of all we observe, thanks to property (ix) of Lemma 3.5, that for $1 \le j \le d$,

$$(4.2)\quad \begin{aligned} &n^{1/2} \left\{ (K^j)^{-1/2} \left( \frac{1}{n} \sum_{k=1}^n \mathbf{S}_{k-1}(j)(\widehat{a}_k^j - a^j)^*(\widehat{a}_k^j - a^j) \right)(K^j)^{-1/2} - \mathbf{I}_d \right\} \\ &= n^{1/2} \left\{ \left( \frac{1}{n} \sum_{k=1}^n (\mathbf{S}_{k-1}(j)^{-1} u(j)) \mathbf{M}_k^{j*} \mathbf{M}_k^j \right) - \mathbf{I}_d \right\} \\ &= n^{1/2} \left( \frac{1}{n} \sum_{k=1}^n (\mathbf{S}_{k-1}(j)^{-1} u(j)) (\mathbf{M}_k^{j*} \mathbf{M}_k^j - \langle \mathbf{M}^j \rangle_k) \right) \\ &= \left( \frac{\mathbf{W}}{\rho - 1} \right)^{-1} n^{1/2} \left\{ \frac{1}{n} \sum_{k=1}^n \rho^{-k} (\mathbf{M}_k^{j*} \mathbf{M}_k^j - \langle \mathbf{M}^j \rangle_k) \right\} + \Delta_n(j), \end{aligned}$$

where $\Delta_n(j) \xrightarrow[n\to\infty]{} 0$ a.s. on $\mathbf{E}$. Hence, according to property (P$'_9$) of Lemma 4.4, we can affirm that conditional on $\mathbf{E}$ or $\mathbf{E}_n$,

$$\left\{ \left[ \sqrt{n} \left( \frac{1}{n} \sum_{k=1}^n \mathbf{S}_{k-1}(j)(\widehat{a}_k^j - a^j)^*(\widehat{a}_k^j - a^j) - K^j \right) \right]_{j\le d}, \mathcal{S}_{n-1}^{1/2}(\widehat{\mathcal{A}}_n - \mathcal{A}) \right\}$$

(4.3)
$$\xrightarrow[n\to\infty]{\mathcal{L}}\left\{\left[\Sigma_1^j\left(\sqrt{\frac{2}{\rho-1}}K^j\right)+\Sigma_2^j(K^j)\right]_{j\le d},\mathfrak{G}\right\},$$

where $\mathfrak{G}$ is a r.v. distributed as $\mathcal{N}_{d^2}(0,\mathcal{K})$, independently of the r.v.'s $\Sigma_r^j(T)$, $1\le j\le d$, $r\in\{1,2\}$, defined as in (P'$_9$).

Now, we notice that we can replace $\frac{1}{n}\sum_{k=1}^n\mathbf{S}_{k-1}(j)(\hat{a}_k^j-a^j)^*(\hat{a}_k^j-a^j)$ by $\widehat{K}_n^j$ in (4.3), because we have

(4.4)
$$\sqrt{n}\left|\frac{1}{n}\sum_{k=1}^n\mathbf{S}_{k-1}(j)(\hat{a}_k^j-a^j)^*(\hat{a}_k^j-a^j)-\widehat{K}_n^j\right|=O\left(\frac{1}{\sqrt{n}}\ln n\right)$$

a.s. on $\mathbf{E}$,

thanks to property (viii) stated in the proof of Theorem 3.2. Hence, conditional on $\mathbf{E}$ or $\mathbf{E}_n$, we get the property

$$\{[\sqrt{n}(\widehat{K}_n^j-K^j)]_{j\le d},\mathcal{S}_{n-1}^{1/2}(\widehat{\mathcal{A}}_n-\mathcal{A})\}$$
$$\xrightarrow[n\to\infty]{\mathcal{L}}\left\{\left[\Sigma_1^j\left(\left(\frac{2}{\rho-1}\right)^{1/4}K^j\right)+\Sigma_2^j(K^j)\right]_{j\le d},\mathfrak{G}\right\},$$

which implies that

$$\{\sqrt{n(\rho-1)}(\widehat{\mathcal{K}}_n-\mathcal{K}),\mathcal{S}_{n-1}^{1/2}(\widehat{\mathcal{A}}_n-\mathcal{A})\}\xrightarrow[n\to\infty]{\mathcal{L}}\{\mathfrak{G}_1(\sqrt{2}\mathcal{K})+\mathfrak{G}_2(\sqrt{\rho-1}\mathcal{K}),\mathfrak{G}\},$$

where, for $\mathcal{T}=\text{Diag}(T^1,\ldots,T^d)$, $\mathfrak{G}_1(\mathcal{T})$ and $\mathfrak{G}_2(\mathcal{T})$ are independent identically distributed Gaussian matrices, which are independent of $\mathfrak{G}$. The covariance matrix of $\text{Vect}\,\mathfrak{G}_r(\mathcal{T})$ is

$$\text{Diag}(T^1\otimes T^1+{}^\perp(\text{Vect}\,T^{1*}\text{Vect}\,T^1),\ldots,T^d\otimes T^d+{}^\perp(\text{Vect}\,T^{d*}\text{Vect}\,T^d)).$$

Property (xi) is proved.

Property (xii) of the theorem, for the sequence $(\widehat{\mathcal{K}}_n)$, is a consequence of property (xv) of Lemma 4.5 combined with (4.4), which implies that

$$\sqrt{n}\left(\text{tr}(\mathcal{K}^{-1/2}\widehat{\mathcal{K}}_n\mathcal{K}^{-1/2})-\frac{1}{n}\sum_{k=1}^n\|\mathcal{K}^{-1/2}\mathcal{S}_{k-1}^{1/2}(\widehat{\mathcal{A}}_k-\mathcal{A})\|^2\right)=O\left(\frac{\ln n}{\sqrt{n}}\right)$$

a.s. on $\mathbf{E}$.

The proofs of properties (xi) and (xii) for the sequence $(\widecheck{\mathcal{K}}_n)$ are similar to those of $(\widehat{\mathcal{K}}_n)$. They are omitted for brevity. $\square$

**5. About the least squares estimator of the mean matrix.** The main contribution of this work is the global identification of the means and the covariance matrices of the reproduction distribution involved in a BGW($d$) process. It was carried out thanks to the CLT of the QSL verified by the normalized estimation errors $(\mathcal{S}_{n-1}(\widehat{\mathcal{A}}_n - \mathcal{A}))$ or $(\mathcal{X}_{n-1}(\widehat{\mathscr{A}}_n - \mathcal{A}))$, where $\widehat{\mathcal{A}}_n$ or $\widehat{\mathscr{A}}_n$ is, respectively, the maximum likelihood or the empirical estimator of the reproduction law means. One may ask if it is possible to adapt this method by considering the least squares estimator (LSE) (see [29]) $\widetilde{\mathbf{A}}_n$ of the mean matrix $\mathbf{A} = [a^1, \ldots, a^d]$ instead of $\widehat{\mathcal{A}}_n$ or $\widehat{\mathscr{A}}_n$. This question is quite relevant if the observed sample is the set of the first $(n+1)$ observations $\mathbf{X}_0, \ldots, \mathbf{X}_n$.

At the beginning, let us discuss the simple case of the BGW process. The LSE of the mean law reproduction $a$ [obtained by minimizing the function $a \longmapsto \sum_{k=1}^n (\mathbf{X}_k - a\mathbf{X}_{k-1})^2]$ is given by

$$(5.1) \qquad \widetilde{a}_n = (\mathbf{Q}_{n-1})^{-1}\left(\sum_{k=1}^n \mathbf{X}_{k-1}\mathbf{X}_k\right), \qquad \mathbf{Q}_n = \sum_{k=0}^n \mathbf{X}_k^2,$$

and satisfies the relation

$$(5.2) \qquad \begin{aligned} \mathbf{Q}_{n-1}(\widetilde{a}_n - a) &= \sum_{k=1}^n \mathbf{X}_{k-1}(\mathbf{X}_k - a\mathbf{X}_{k-1}) \\ &= \sum_{k=1}^n \mathbf{X}_{k-1}(\mathbf{X}_k - \mathbb{E}(\mathbf{X}_k/\mathcal{F}_{k-1})) = \mathbf{W}_n. \end{aligned}$$

Noting that the predictable quadratic variation of the martingale $(\mathbf{W}_n)$ is

$$\langle \mathbf{W} \rangle_n = \sigma^2 \mathbf{T}_{n-1} \qquad \text{with } \mathbf{T}_n = \sum_{k=0}^n \mathbf{X}_k^3,$$

we can affirm, thanks to the martingale law of large numbers, that $(\widetilde{a}_n)$ is a strongly consistent estimator of $a$ on the set $\mathbf{E}$. Moreover, conditional on the set $\{\mathbf{X}_n > 0\}$ (resp. $\mathbf{E}$), the following CLT with random normalization holds too:

$$(5.3) \qquad \begin{aligned} \sqrt{\mathbf{S}_{n-1}}(\widetilde{a}_n - a) &\xrightarrow[n \to \infty]{\mathcal{L}} \mathcal{N}(0, \gamma^2), \\ \sqrt{\mathbf{T}_{n-1}}\mathbf{Q}_{n-1}(\widetilde{a}_n - a) &\xrightarrow[n \to \infty]{\mathcal{L}} \mathcal{N}(0, \sigma^2), \end{aligned}$$

where $\gamma^2 = ((a+1)^2/(a^2+a+1))\sigma^2 > \sigma^2$.

The QSL associated with this CLT allows us to prove that

$$(5.4) \qquad \widetilde{\sigma}_n^2 = n^{-1}\sum_{k=1}^n \mathbf{T}_{k-1}\mathbf{Q}_{k-1}^2(\widetilde{a}_k - \widetilde{a}_n)^2 \xrightarrow[n \to \infty]{} \sigma^2 \qquad \text{a.s. on } \mathbf{E}.$$

Using a CLT stated in [8], we have an analogous result to (1.5), which is

$$(5.5) \left\{ \sqrt{\frac{n(a^3-1)}{2(a^3+1)}} \left( \left(\frac{\widetilde{\sigma}_n}{\sigma}\right)^2 - 1 \right), \sqrt{\frac{\mathbf{S}_{n-1}}{\widetilde{\sigma}_n^2}} (\widehat{a}_n - a) \right\} \xrightarrow[n\to\infty]{\mathcal{L}} \mathcal{N}(0,1) \otimes \mathcal{N}(0,1).$$

Obviously, the estimators $(\widehat{a}_n, \widetilde{\sigma}_n^2)$ are better than $(\widetilde{a}_n, \widehat{\sigma}_n^2)$. However, as stated before, for the BGW($d$) process, the maximum likelihood or the empirical estimator of the reproduction mean matrix $\mathbf{A}$ is not a function of the first $(n+1)$ observations of the process, but is a function of a much richer sample drawn from the underlying tree generated by this process. So it is interesting to study the properties of the LSE $(\widetilde{A}_n)$ of $\mathbf{A}$ because it is a statistic of the basic sample $(\mathbf{X}_0, \ldots, \mathbf{X}_n)$. More precisely, $\widetilde{A}_n$ (calculated by minimizing the function $A \longmapsto \sum_{k=1}^n \|\mathbf{X}_k - A\mathbf{X}_{k-1}\|^2$) is a solution of the linear system

$$\widetilde{A}_n \left( \sum_{k=1}^n \mathbf{X}_{k-1}^* \mathbf{X}_{k-1} \right) = \sum_{k=1}^n \mathbf{X}_k^* \mathbf{X}_{k-1}.$$

It is more convenient to define $\widetilde{A}_n$ by the algorithm

$$(5.6)\, \widetilde{A}_{n+1} = \widetilde{A}_n + (\mathbf{X}_{n+1} - \widetilde{A}_n \mathbf{X}_n)^* \mathbf{X}_n \mathbf{Q}_n^{-1}, \qquad \text{where } \mathbf{Q}_n = \mathbf{I}_d + \sum_{k=1}^n \mathbf{X}_k^* \mathbf{X}_k.$$

Hence, $\widetilde{A}_n$ satisfies the relation analogous to (5.2),

$$(5.7) \qquad (\widetilde{A}_n - A)\mathbf{Q}_{n-1} = \sum_{j=1}^d \mathbf{W}_n^j, \qquad \text{where } \mathbf{W}_n^j = \sum_{k=1}^n \mathbf{X}_{k-1}^* \zeta_k^j.$$

The global approach discussed before for the one-type BGW can be generalized successfully to the $d$-type case. For this purpose, we need to add to Assumptions A-1 and A-2:

ASSUMPTION A-3. Every nonprincipal eigenvalue $\lambda$ of $\mathbf{A}$ satisfies $|\lambda|^2 > \rho$. Moreover, $\mathbf{A}$ is nonderogatory.

The last word means that the minimal and characteristic polynomials of $\mathbf{A}$ are proportional.

Under Assumption A-2 and the first part of Assumption A-3, it was proved by Carvalho [5] that there exists a r.v. $\eta$ that satisfies the property

$$(5.8) \qquad \mathbf{A}^{-n}\mathbf{X}_n \xrightarrow[n\to\infty]{} \eta \qquad \text{a.s. (and also in mean square).}$$

Moreover, when we add the second part of Assumption A-3, we get the important property (stated also in [5]):

$$(5.9) \qquad \begin{array}{l} \text{On the set of nonextinction } \mathbf{E}, \text{ the } d \times d \text{ matrix whose} \\ \text{column vectors are } [\mathbf{A}^{-1}\eta, \ldots, \mathbf{A}^{-d}\eta] \text{ is invertible.} \end{array}$$

As a direct consequence of (5.9), the excitation $\mathbf{Q}_n$ of $\mathbf{X}$ satisfies the result

$$\mathbf{A}^{-n}\mathbf{Q}_{n-1}{}^*\mathbf{A}^{-n} \underset{n\to\infty}{\longrightarrow} \mathrm{C} \qquad \text{a.s. on } \mathbf{E},$$

where $\mathbf{C} = \sum_{n=1}^{\infty} \mathbf{A}^{-n}\eta^*\eta^*\mathbf{A}^{-n}$ is a.s. invertible on $\mathbf{E}$.

By adapting to the martingale

$$\widetilde{\mathcal{M}}_n = {}^*(\mathrm{Vect}(\mathbf{W}_n^1), \dots, \mathrm{Vect}(\mathbf{W}_n^d))$$

the tools used before for $(\mathcal{M}_n)$, similar results to those proved for the MLE $\widehat{\mathcal{A}}_n$ (or the EE $\widehat{\mathcal{A}}_n$) are available for $\widetilde{\mathcal{A}}_n = \mathrm{Vect}(\widetilde{A}_n)$. Details are intentionally omitted for brevity.

**Conclusion.** A natural question remains: What is the best estimator of the variance? The answer to this question is part of a thesis in preparation.

## APPENDIX A: PROOFS OF LEMMAS 3.4 AND 3.5

PROOF OF LEMMA 3.4. The ASCLT satisfied by $(\mathcal{M}_n)$ or $(\zeta_n)$ has been proved in [10]. Property $(P_7)$ for the martingale $(\mathcal{M}_n)$ is a special case of Theorem 3.1 in [8]. Property $(P_8)$ for the martingale $(\mathcal{M}_n)$ is proved at the end of Appendix B. The proofs of the properties $(P_7)$ and $(P_8)$ for the sequence $(\zeta_n)$ are similar to those for $(\mathcal{M}_n)$. They are omitted for brevity. □

PROOF OF LEMMA 3.5. To prove property (ix) of the lemma, we note that

$$n^{-1/2} \sum_{k=1}^{n} (\mathcal{K}^{-1/2}\mathcal{S}_{k-1}^{1/2}(\widehat{\mathcal{A}}_k - \mathcal{A})^*(\widehat{\mathcal{A}}_k - \mathcal{A})\mathcal{S}_{k-1}^{1/2}\mathcal{K}^{-1/2} - \mathcal{I}_{d^2})$$

$$\begin{aligned}
\text{(A.1)} \quad &= n^{-1/2} \sum_{k=1}^{n} (\rho^k \mathcal{U}\mathcal{S}_{k-1}^{-1})^{1/2} \rho^{-k}(\mathcal{M}_k{}^*\mathcal{M}_k - \langle\mathcal{M}\rangle_k)(\rho^k \mathcal{U}\mathcal{S}_{k-1}^{-1})^{1/2} \\
&= \left(\frac{\mathbf{W}}{\rho-1}\right)^{-1} n^{-1/2} \sum_{k=1}^{n} \rho^{-k}(\mathcal{M}_k{}^*\mathcal{M}_k - \langle\mathcal{M}\rangle_k) + \Gamma_n,
\end{aligned}$$

where

$$\begin{aligned}
\Gamma_n = n^{-1/2} \sum_{k=1}^{n} &\left((\rho^k \mathcal{U}\mathcal{S}_{k-1}^{-1})^{1/2} - \left(\frac{\mathbf{W}}{\rho-1}\right)^{-1/2}\mathcal{I}_{d^2}\right) \\
&\times \rho^{-k}(\mathcal{M}_k{}^*\mathcal{M}_k - \langle\mathcal{M}\rangle_k)\left((\rho^k \mathcal{U}\mathcal{S}_{k-1}^{-1})^{1/2} - \left(\frac{\mathbf{W}}{\rho-1}\right)^{-1/2}\mathcal{I}_{d^2}\right)
\end{aligned}$$

$$+ \left(\frac{\mathbf{W}}{\rho-1}\right)^{-1/2} n^{-1/2} \sum_{k=1}^{n} \left((\rho^k \mathcal{U} \mathcal{S}_{k-1}^{-1})^{1/2} - \left(\frac{\mathbf{W}}{\rho-1}\right)^{-1/2} \mathcal{I}_{d^2}\right)$$

$$\times \rho^{-k}(\mathcal{M}_k{}^*\mathcal{M}_k - \langle\mathcal{M}\rangle_k)$$

$$+ \left(\frac{\mathbf{W}}{\rho-1}\right)^{-1/2} n^{-1/2} \sum_{k=1}^{n} \rho^{-k}(\mathcal{M}_k{}^*\mathcal{M}_k - \langle\mathcal{M}\rangle_k)$$

$$\times \left((\rho^k \mathcal{U} \mathcal{S}_{k-1}^{-1})^{1/2} - \left(\frac{\mathbf{W}}{\rho-1}\right)^{-1/2} \mathcal{I}_{d^2}\right).$$

Let us show that almost surely on $\mathbf{E}$,

(A.2) $$\Gamma_n \underset{n\to\infty}{\longrightarrow} 0.$$

Properties $(P_2)$ and $(P_3)$ imply $\rho^{-n}\|\mathcal{M}_n{}^*\mathcal{M}_n - \langle\mathcal{M}\rangle_n\| = O(\ln n) = o(n)$ a.s. on $\mathbf{E}$, so

$$\Gamma_n = O\left(n^{-1/2} \sum_{k=1}^{n} k \left\|(\rho^k \mathcal{U} \mathcal{S}_{k-1}^{-1})^{1/2} - \left(\frac{\mathbf{W}}{\rho-1}\right)^{-1/2} \mathcal{I}_{d^2}\right\|^2\right)$$

(A.3) $$+ O\left(n^{-1/2} \sum_{k=1}^{n} k \left\|(\rho^k \mathcal{U} \mathcal{S}_{k-1}^{-1})^{1/2} - \left(\frac{\mathbf{W}}{\rho-1}\right)^{-1/2} \mathcal{I}_{d^2}\right\|\right)$$

a.s. on $\mathbf{E}$.

To end the proof of (A.2), we need the following result, which is a direct consequence of property $(P'_1)$: For some real $\theta \in ]0,1[$,

(A.4) $$\left\|\rho^{-k}\mathcal{S}_{k-1} - \frac{1}{\rho-1}\mathbf{W}\mathcal{U}\right\| = O(\sqrt{\ln k}\, k^{d-1}\theta^k) \qquad \text{a.s. on } \mathbf{E}.$$

This implies that almost surely on $\mathbf{E}$,

$$\left\|(\rho^k \mathcal{U} \mathcal{S}_{k-1}^{-1})^{1/2} - \left(\frac{\mathbf{W}}{\rho-1}\right)^{-1/2} \mathcal{I}_{d^2}\right\| = O(\sqrt{\ln k}\, k^{d-1}\theta^k).$$

Therefore, the series

$$\sum_{k=1}^{\infty} \sqrt{k} \left\|(\rho^k \mathcal{U} \mathcal{S}_{k-1}^{-1})^{1/2} - \left(\frac{\mathbf{W}}{\rho-1}\right)^{-1/2} \mathcal{I}_{d^2}\right\|^2$$

and

$$\sum_{k=1}^{\infty} \sqrt{k} \left\|(\rho^k \mathcal{U} \mathcal{S}_{k-1}^{-1})^{1/2} - \left(\frac{\mathbf{W}}{\rho-1}\right)^{-1/2} \mathcal{I}_{d^2}\right\|$$

are a.s. convergent on $\mathbf{E}$. By the Kronecker lemma, almost surely on $\mathbf{E}$,

$$\sum_{k=1}^{n} k \left\| (\rho^k \mathcal{U} \mathcal{S}_{k-1}^{-1})^{1/2} - \left( \frac{\mathbf{W}}{\rho-1} \right)^{-1/2} \mathcal{I}_{d^2} \right\|^2$$

$$+ \sum_{k=1}^{n} k \left\| (\rho^k \mathcal{U} \mathcal{S}_{k-1}^{-1})^{1/2} - \left( \frac{\mathbf{W}}{\rho-1} \right)^{-1/2} \mathcal{I}_{d^2} \right\| = o(\sqrt{n});$$

hence (A.2) is proved thanks to (A.3).

Since property (x) is an immediate consequence of (ix), Lemma 3.5 is proved for $(\mathcal{S}_{n-1}^{1/2}(\widehat{\mathcal{A}}_n - \mathcal{A}))$.

The proofs of properties (ix) and (x) for $(\mathcal{X}_{n-1}^{1/2}(\widehat{\widetilde{\mathcal{A}}}_n - \mathcal{A}))$ are similar to those for $(\mathcal{S}_{n-1}^{1/2}(\widehat{\mathcal{A}}_n - \mathcal{A}))$. They are omitted for brevity.  $\square$

## APPENDIX B: PROOFS OF LEMMA 4.4 AND PROPERTY (P$_8$)

PROOF OF LEMMA 4.4.   The long and technical proof below is carried in several steps. These steps are intentionally detailed in a way so that results on the sequence $(\zeta_n)$ appear to be contained in those on the martingale $(\mathcal{M}_n)$; hence the proof will be focused on $(\mathcal{M}_n)$.

STEP 1—Preliminary calculus.   Let

$$\mathcal{Z}_n = \sum_{k=1}^{n} \rho^{-k} (\mathcal{M}_k{}^* \mathcal{M}_k - \langle \mathcal{M} \rangle_k),$$

(B.1)

$$\widetilde{\mathcal{Z}}_n = \sum_{k=1}^{n} \mathcal{K}^{1/2} \mathcal{X}_{k-1}^{-1/2} (\zeta_k{}^* \zeta_k - \mathcal{X}_{k-1}^{1/2} \mathcal{K} \mathcal{X}_{k-1}^{1/2}) \mathcal{X}_{k-1}^{-1/2} \mathcal{K}^{-1/2}$$

be the random block matrices whose blocks of indexes $i$ and $r$, $1 \le i, r \le d$, are

$$\mathcal{Z}_n(i,r) = \sum_{k=1}^{n} \rho^{-k} (\mathbf{M}_k^{i}{}^* \mathbf{M}_k^r - \langle \mathbf{M}^i, \mathbf{M}^r \rangle_k),$$

$$\widetilde{\mathcal{Z}}_n(i,r) = \sum_{k=1}^{n} (\mathbf{X}_{k-1}(i)^{-1/2} \mathbf{X}_{k-1}(r)^{-1/2} (K^i)^{-1/2} \zeta_k^{i*} \zeta_k^r (K^r)^{-1/2} - \mathbf{I}_d).$$

To prove property (P$_9$), we need the decomposition

$$\mathcal{Z}_n(i,r) = \left( \frac{\rho}{\rho-1} \right) \sum_{k=1}^{n} \left( \frac{1}{\rho^k} - \frac{1}{\rho^{k+1}} \right) (\mathbf{M}_k^{i*} \mathbf{M}_k^r - \langle \mathbf{M}^i, \mathbf{M}^r \rangle_k)$$

$$= \left( \frac{\rho}{\rho-1} \right) \sum_{k=1}^{n} \rho^{-k} (\mathbf{M}_{k-1}^{i*} (\Delta \mathbf{M}_k^r) + (\Delta \mathbf{M}_k^i)^* \mathbf{M}_{k-1}^r)$$

$$+ \left(\frac{\rho}{\rho-1}\right) \sum_{k=1}^{n} \rho^{-k}((\Delta \mathbf{M}_k^i)^*(\Delta \mathbf{M}_k^r) - \mathbb{E}((\Delta \mathbf{M}_k^i)^*(\Delta \mathbf{M}_k^r)|\mathcal{G}_{k-1}))$$

$$- \left(\frac{\rho}{\rho-1}\right) \rho^{-(n+1)}(\mathbf{M}_n^{i*}\mathbf{M}_n^r - \langle \mathbf{M}^i, \mathbf{M}^r \rangle_n).$$

From $(\mathrm{P}_2)$ and $(\mathrm{P}_3)$, we have

$$n^{-1/2}\left(\frac{\rho}{\rho-1}\right) \rho^{-(n+1)}(\mathcal{M}_n^*\mathcal{M}_n - \langle \mathcal{M} \rangle_n) \underset{n \to \infty}{\longrightarrow} 0 \qquad \text{a.s. on } \mathbf{E};$$

hence,

$$(\mathrm{B.2}) \qquad \left(\frac{\rho-1}{\rho}\right) \mathcal{Z}_n(i,r) = H_n^{i,r} + \widetilde{H}_n^{i,r} + o(\sqrt{n}) \qquad \text{a.s. on } \mathbf{E},$$

where for $1 \le i, r \le d$,

$$H_n^{i,r} = \sum_{k=1}^{n} \rho^{-k}(\mathbf{M}_{k-1}^{i*}(\Delta \mathbf{M}_k^r) + (\Delta \mathbf{M}_k^i)^*\mathbf{M}_{k-1}^r),$$

$$\widetilde{H}_n^{i,r} = \sum_{k=1}^{n} \rho^{-k}((\Delta \mathbf{M}_k^i)^*(\Delta \mathbf{M}_k^r) - \mathbb{E}((\Delta \mathbf{M}_k^i)^*(\Delta \mathbf{M}_k^r)|\mathcal{G}_{k-1}))$$

$$= \frac{\mathbf{W}}{\rho} \widetilde{\mathcal{Z}}_n(i,r) + O(\ln n) \qquad \text{a.s. on } \mathbf{E}.$$

For the study of the weak convergence of the sequence $(\mathcal{Z}_n)$, we consider the martingale block matrices

$$\mathsf{H}_n = (H_n^{i,r})_{1 \le i, r \le d}, \qquad \widetilde{\mathsf{H}}_n = (\widetilde{H}_n^{i,r})_{1 \le i, r \le d},$$

$$\mathcal{H}_n = \mathrm{Vect}(\mathsf{H}_n) = \mathrm{Vect}(\mathsf{H}_n^1, \ldots, \mathsf{H}_n^{d^2}),$$

$$\widetilde{\mathcal{H}}_n = \mathrm{Vect}(\widetilde{\mathsf{H}}_n) = \mathrm{Vect}(\widetilde{\mathsf{H}}_n^1, \ldots, \widetilde{\mathsf{H}}_n^{d^2}) \quad \text{and} \quad \mathbb{H}_n = \begin{pmatrix} \mathcal{H}_n \\ \widetilde{\mathcal{H}}_n \end{pmatrix}.$$

The continuation of the proof is, in particular, based on application of the classical CLT (see [15]) to the martingale $(\mathbb{H}_n)$.

STEP 2—Behavior of the predictable quadratic variation of $(\mathbb{H}_n)$.  Hereafter, we study the asymptotic behavior of the predictable quadratic variation

$$\langle \mathbb{H} \rangle_n = \begin{pmatrix} \langle \mathcal{H} \rangle_n & \langle \mathcal{H}, \widetilde{\mathcal{H}} \rangle_n \\ {}^*\langle \mathcal{H}, \widetilde{\mathcal{H}} \rangle_n & \langle \widetilde{\mathcal{H}} \rangle_n \end{pmatrix} \text{ of } \mathbb{H}_n,$$

where

$$\langle \mathcal{H} \rangle_n = (\langle \mathsf{H}^p, \mathsf{H}^q \rangle_n)_{p,q \le d^2},$$

$$\langle \widetilde{\mathcal{H}} \rangle_n = (\langle \widetilde{\mathsf{H}}^p, \widetilde{\mathsf{H}}^q \rangle_n)_{p,q \le d^2},$$

$$\langle \mathcal{H}, \widetilde{\mathcal{H}} \rangle_n = (\langle \mathsf{H}^p, \widetilde{\mathsf{H}}^q \rangle_n)_{p,q \le d^2}.$$

To determine the asymptotic behavior of $(\langle \mathbb{H} \rangle_n)$, the formulas below are needed, where $(e_1, \ldots, e_d)$ denotes the canonical basis of $\mathbb{R}^d$. For $1 \leq r \leq d$, $1 \leq s \leq d$, we note that

$$\mathsf{H}_n^{(r-1)d+s} = {}^*({}^*H_n^{1,r}e_s, \ldots, {}^*H_n^{d,r}e_s),$$

$$\widetilde{\mathsf{H}}_n^{(r-1)d+s} = {}^*({}^*\widetilde{H}_n^{1,r}e_s, \ldots, {}^*\widetilde{H}_n^{d,r}e_s).$$

Then, for $p = (r-1)d + s$ and $q = (l-1)d + t$ with $1 \leq r, l \leq d$ and $1 \leq s, t \leq d$ fixed,

$$\langle \mathsf{H}^p, \mathsf{H}^q \rangle_n = (\langle H_n^{i,r}e_s, H_n^{j,l}e_t \rangle_n)_{1 \leq i,j \leq d},$$

$$\langle \mathsf{H}^p, \widetilde{\mathsf{H}}^q \rangle_n = (\langle H_n^{i,r}e_s, \widetilde{H}_n^{j,l}e_t \rangle_n)_{1 \leq i,j \leq d},$$

$$\langle \widetilde{\mathsf{H}}^p, \widetilde{\mathsf{H}}^q \rangle_n = (\langle \widetilde{H}_n^{i,r}e_s, \widetilde{H}_n^{j,l}e_t \rangle_n)_{1 \leq i,j \leq d}.$$

(A) *Behavior of* $\langle \mathcal{H} \rangle_n$. For $1 \leq i, j \leq d$, $1 \leq r, l \leq d$, $1 \leq s, t \leq d$, we have $\langle H^{i,r}e_s, H^{j,l}e_t \rangle_n$

$$\begin{aligned}
&= \sum_{k=1}^n \rho^{-2k} \mathbf{M}_{k-1}^i {}^* \mathbf{M}_{k-1}^j \mathbb{E}(\langle \mathbf{\Delta M}_k^r, e_s \rangle \langle \mathbf{\Delta M}_k^l, e_t \rangle | \mathcal{G}_{k-1}) \\
&\quad + \sum_{k=1}^n \rho^{-2k} \mathbf{M}_{k-1}^i \langle \mathbf{M}_{k-1}^l, e_t \rangle \mathbb{E}(\langle \mathbf{\Delta M}_k^r, e_s \rangle^* \mathbf{\Delta M}_k^j | \mathcal{G}_{k-1}) \\
&\quad + \sum_{k=1}^n \rho^{-2k} \mathbb{E}(\mathbf{\Delta M}_k^i \langle \mathbf{\Delta M}_k^l, e_t \rangle | \mathcal{G}_{k-1}) \langle \mathbf{M}_{k-1}^r, e_s \rangle^* \mathbf{M}_{k-1}^j \\
&\quad + \sum_{k=1}^n \rho^{-2k} \langle \mathbf{M}_{k-1}^r, e_s \rangle \langle \mathbf{M}_{k-1}^l, e_t \rangle \mathbb{E}((\mathbf{\Delta M}_k^i)^*(\mathbf{\Delta M}_k^j) | \mathcal{G}_{k-1}),
\end{aligned}$$

$\langle H^{i,r}e_s, H^{j,l}e_t \rangle_n$

$$\begin{aligned}
&= \frac{1}{\sqrt{u(r)u(l)}} \sum_{k=1}^n \rho^{-2k} \mathbf{M}_{k-1}^i {}^* \mathbf{M}_{k-1}^j {}^* e_s (K^r)^{-1/2} \mathbb{E}(\zeta_k^{r*} \zeta_k^l | \mathcal{G}_{k-1})(K^l)^{-1/2} e_t \\
&\quad + \frac{1}{\sqrt{u(r)u(j)}} \sum_{k=1}^n \rho^{-2k} \mathbf{M}_{k-1}^i {}^* \mathbf{M}_{k-1}^l e_t {}^* e_s (K^r)^{-1/2} \\
&\qquad\qquad \times \mathbb{E}(\zeta_k^{r*} \zeta_k^j | \mathcal{G}_{k-1})(K^j)^{-1/2} \\
&\quad + \frac{1}{\sqrt{u(i)u(l)}} \sum_{k=1}^n \rho^{-2k} (K^i)^{-1/2} \mathbb{E}(\zeta_k^{i*} \zeta_k^l | \mathcal{G}_{k-1})(K^l)^{-1/2} e_t {}^* e_s \mathbf{M}_{k-1}^r {}^* \mathbf{M}_{k-1}^j \\
&\quad + \frac{1}{\sqrt{u(i)u(j)}} \sum_{k=1}^n \rho^{-2k} {}^* e_s \mathbf{M}_{k-1}^r {}^* \mathbf{M}_{k-1}^l e_t (K^i)^{-1/2} \\
&\qquad\qquad \times \mathbb{E}(\zeta_k^{i*} \zeta_k^j | \mathcal{G}_{k-1})(K^j)^{-1/2}.
\end{aligned}$$

Consequently,

$$
\begin{aligned}
\langle H^{i,r}e_s, H^{j,l}e_t\rangle_n = {}& \langle e_r, e_l\rangle^* e_s e_t\left(\sum_{k=1}^n u(r)^{-1}\mathbf{X}_{k-1}(r)\rho^{-2k}\mathbf{M}_{k-1}^{i}{}^*\mathbf{M}_{k-1}^{j}\right) \\
& + \langle e_j, e_r\rangle\left(\sum_{k=1}^n u(r)^{-1}\mathbf{X}_{k-1}(r)\rho^{-2k}\mathbf{M}_{k-1}^{i}{}^*\mathbf{M}_{k-1}^{l}\right)e_t{}^*e_s \\
& + \langle e_i, e_l\rangle\left(\sum_{k=1}^n u(l)^{-1}\mathbf{X}_{k-1}(l)\rho^{-2k}e_t{}^*e_s\mathbf{M}_{k-1}^{r}{}^*\mathbf{M}_{k-1}^{j}\right) \\
& + \langle e_i, e_j\rangle\left(\sum_{k=1}^n u(i)^{-1}\mathbf{X}_{k-1}(i)\rho^{-2k}{}^*e_s\mathbf{M}_{k-1}^{r}{}^*\mathbf{M}_{k-1}^{l}e_t\mathbf{I}_d\right).
\end{aligned}
$$
(B.3)

From (B.3) and the properties

$$
\rho^{-(n-1)}\mathbf{X}_{n-1}(r)\underset{n\to\infty}{\longrightarrow}\mathbf{W}u(r),
$$

(B.4)
$$
\frac{1}{n}\sum_{k=1}^n \rho^{-(k-1)}\mathbf{M}_{k-1}^{r}{}^*\mathbf{M}_{k-1}^{l}\underset{n\to\infty}{\longrightarrow}\left(\frac{\mathbf{W}}{\rho-1}\mathbf{I}_d\right)\langle e_r, e_l\rangle \qquad \text{a.s. on } \mathbf{E},
$$

we can affirm that, almost surely on $\mathbf{E}$,

$$
\frac{1}{n}\langle H^{i,r}e_s, H^{j,l}e_t\rangle_n \longrightarrow \frac{2\mathbf{W}^2}{\rho^2(\rho-1)}({}^*e_i\mathbf{I}_d e_j\langle e_r, e_l\rangle\langle e_s, e_t\rangle\mathbf{I}_d + {}^*e_i(e_l{}^*e_r)e_j(e_t{}^*e_s)).
$$

So, we deduce that

$$
\frac{1}{n}\langle \mathsf{H}^p, \mathsf{H}^q\rangle_n \underset{n\to\infty}{\longrightarrow}\frac{2\mathbf{W}^2}{\rho^2(\rho-1)}(\langle e_r, e_l\rangle\langle e_s, e_t\rangle\mathbf{I}_d\otimes\mathbf{I}_d + (e_l{}^*e_r)\otimes(e_s{}^*e_s))
$$
$$
\text{a.s. on } \mathbf{E}
$$

and

$$
\frac{1}{n}\langle \mathcal{H}\rangle_n\underset{n\to\infty}{\longrightarrow}\frac{2\mathbf{W}^2}{\rho^2(\rho-1)}\{(\mathbf{I}_d\otimes\mathbf{I}_d)\otimes(\mathbf{I}_d\otimes\mathbf{I}_d) + \mathcal{J}\boxtimes\mathcal{J}\} \qquad \text{a.s. on } \mathbf{E},
$$

where

(B.5)
$$
\mathcal{J} = {}^\perp(\mathrm{Vect}(\mathbf{I}_d)^*\,\mathrm{Vect}(\mathbf{I}_d)) = (e_l{}^*e_r)_{1\leq r,l\leq d}.
$$

(B) *Behavior of* $(\langle\widetilde{\mathcal{H}}\rangle_n)$. For $1\leq i,j\leq d, 1\leq r,l\leq d, 1\leq s,t\leq d$, we have

$$
\begin{aligned}
\langle \widetilde{H}^{i,r}e_s, \widetilde{H}^{j,l}e_t\rangle_n \\
= \sum_{k=1}^n \rho^{-2k}\mathbb{E}((\Delta\mathbf{M}_k^{i})^*(\Delta\mathbf{M}_k^{j})\langle\Delta\mathbf{M}_k^{r}, e_s\rangle\langle\Delta\mathbf{M}_k^{l}, e_t\rangle|\mathcal{G}_{k-1})
\end{aligned}
$$

$$-\sum_{k=1}^{n}\rho^{-2k}\mathbb{E}((\Delta\mathbf{M}_k^i)\langle\Delta\mathbf{M}_k^r,e_s\rangle|\mathcal{G}_{k-1})\mathbb{E}(^*(\Delta\mathbf{M}_k^j)\langle\Delta\mathbf{M}_k^l,e_t\rangle|\mathcal{G}_{k-1}),$$

$$\langle\widetilde{H}^{i,r}e_s,\widetilde{H}^{j,l}e_t\rangle_n$$

$$=\frac{1}{\sqrt{u(i)u(r)u(j)u(l)}}$$

$$\times\sum_{k=1}^{n}\rho^{-2k}\mathbb{E}((K^i)^{-1/2}\zeta_k^{i*}\zeta_k^j(K^j)^{-1/2}$$

$$\times{}^*e_s(K^r)^{-1/2}\zeta_k^{r*}\zeta_k^l(K^l)^{-1/2}e_t|\mathcal{G}_{k-1})$$

(B.6)     $$-\frac{1}{\sqrt{u(i)u(r)u(j)u(l)}}$$

$$\times\sum_{k=1}^{n}\rho^{-2k}(K^i)^{-1/2}$$

$$\times\mathbb{E}(\zeta_k^{i*}\zeta_k^r|\mathcal{G}_{k-1})(K^r)^{-1/2}e_s{}^*e_t(K^l)^{-1/2}\mathbb{E}(\zeta_k^{l*}\zeta_k^j|\mathcal{G}_{k-1})(K^j)^{-1/2}.$$

Consequently, if we set $\widehat{\xi}_{1,1}^i=(K^i)^{-1/2}(\xi_{1,1}^i-a^j)$, then on the set $\{\mathbf{X}_{k-1}(i)\times\mathbf{X}_{k-1}(j)\mathbf{X}_{k-1}(r)\mathbf{X}_{k-1}(l)\neq 0\}$,

$$\mathbb{E}((K^i)^{-1/2}\zeta_k^{i*}\zeta_k^j(K^j)^{-1/2*}e_s(K^r)^{-1/2}\zeta_k^{r*}\zeta_k^l(K^l)^{-1/2}e_t|\mathcal{G}_{k-1})$$

$$=\mathbb{E}\Bigg(\bigg(\sum_{i_1=1}^{\mathbf{X}_{k-1}(i)}\widehat{\xi}_{k,i_1}^i\bigg)^*\bigg(\sum_{j_1=1}^{\mathbf{X}_{k-1}(j)}\widehat{\xi}_{k,j_1}^j\bigg)^*e_s\bigg(\sum_{r_1=1}^{\mathbf{X}_{k-1}(r)}\widehat{\xi}_{k,r_1}^r\bigg)^*\bigg(\sum_{l_1=1}^{\mathbf{X}_{k-1}(l)}\widehat{\xi}_{k,l_1}^l\bigg)$$

$$\times e_t|\mathcal{G}_{k-1}\Bigg)$$

$$=\sum_{i_1=1}^{\mathbf{X}_{k-1}(i)}\sum_{j_1=1}^{\mathbf{X}_{k-1}(j)}\sum_{r_1=1}^{\mathbf{X}_{k-1}(r)}\sum_{l_1=1}^{\mathbf{X}_{k-1}(l)}\Theta_{(i_1,r_1,j_1,l_1)}^{(i,r,j,l)},$$

where

$$\Theta_{(i_1,r_1,j_1,l_1)}^{(i,r,j,l)}=\mathbb{E}(\widehat{\xi}_{1,i_1}^i{}^*\widehat{\xi}_{1,j_1}^j{}^*e_s\widehat{\xi}_{1,r_1}^r{}^*\widehat{\xi}_{1,l_1}^le_t).$$

To calculate $\Theta_{(i_1,r_1,j_1,l_1)}^{(i,r,j,l)}$ the following three cases must be examined:

*Case* 1.  The four indexes $i$, $r$, $j$, $l$ are equal.
*Case* 2.  Exactly three indexes from $\{i,r,j,l\}$ are equal.
*Case* 3.  Two indexes at most from $\{i,r,j,l\}$ are equal.

In Case 2, $\Theta_{(i_1,r_1,j_1,l_1)}^{(i,r,j,l)}=0$, since for $x\in\mathbb{R}^d,y\in\mathbb{R}^d$,

$$^*x\Theta_{(i_1,r_1,j_1,l_1)}^{(i,i,i,l)}y=\mathbb{E}(\langle x,\widehat{\xi}_{1,i_1}^i\rangle\langle y,\widehat{\xi}_{1,j_1}^i\rangle\langle e_s,\widehat{\xi}_{1,r_1}^i\rangle)\mathbb{E}(\langle\widehat{\xi}_{1,l_1}^l,e_t\rangle)=0.$$

Likewise,

$$^*x\Theta^{(i,i,j,i)}_{(i_1,j_1,r_1,l_1)}y = {}^*x\Theta^{(i,r,i,i)}_{(i_1,j_1,r_1,l_1)}y = {}^*x\Theta^{(i,r,r,r)}_{(i_1,j_1,r_1,l_1)}y = 0.$$

In Case 3, $\Theta^{(i,r,j,l)}_{(i_1,r_1,j_1,l_1)} = 0$, except perhaps in the following subcases:

(a) If $(i,r,j,l) \in \Lambda_1 = \{(r,r,l,l); 1 \le r,l \le d, r \ne l\}$, then

$$\Theta^{(r,r,l,l)}_{(i_1,r_1,j_1,l_1)} = \mathbb{E}(\widehat{\xi}^r_{1,i_1} {}^*\widehat{\xi}^r_{1,r_1})e_s{}^*e_t \mathbb{E}(\widehat{\xi}^l_{1,j_1} {}^*\widehat{\xi}^l_{1,l_1}) = (e_s{}^*e_t)\mathbf{1}_{\Lambda_1}(i_1,r_1,j_1,l_1).$$

(b) If $(i,r,j,l) \in \Lambda_2 = \{(i,r,i,r); 1 \le i,r \le d, i \ne r\}$, then

$$\Theta^{(i,r,i,r)}_{(i_1,r_1,j_1,l_1)} = \mathbb{E}(\widehat{\xi}^i_{1,i_1} {}^*\widehat{\xi}^i_{1,j_1}){}^*e_s \mathbb{E}(\widehat{\xi}^r_{1,r_1} {}^*\widehat{\xi}^r_{1,l_1})e_t = \langle e_s, e_t\rangle \mathbf{I}_d \mathbf{1}_{\Lambda_2}(i_1,r_1,j_1,l_1).$$

(c) If $(i,r,j,l) \in \Lambda_3 = \{(l,r,r,l); 1 \le l,r \le d, l \ne r\}$, then

$$\Theta^{(l,r,r,l)}_{(i_1,r_1,j_1,l_1)} = \mathbb{E}(\widehat{\xi}^l_{1,i_1} {}^*\widehat{\xi}^l_{1,l_1})e_t{}^*e_s \mathbb{E}(\widehat{\xi}^r_{1,r_1} {}^*\widehat{\xi}^r_{1,j_1}) = (e_t{}^*e_s)\mathbf{1}_{\Lambda_3}(i_1,r_1,j_1,l_1).$$

In Case 1, $(i,j,r,l) \in \Lambda_4 = \{(r,r,r,r); 1 \le r \le d\}$. Then

$$\begin{aligned}
\Theta^{(i,j,r,l)}_{(i_1,r_1,j_1,l_1)} &= \mathbb{E}(\widehat{\xi}^r_{1,i_1} {}^*\widehat{\xi}^r_{1,j_1}\langle e_s, \widehat{\xi}^r_{1,r_1}\rangle\langle\widehat{\xi}^r_{1,l_1}, e_t\rangle)\\
&= \mathbb{E}(\widehat{\xi}^r_{1,1} {}^*\widehat{\xi}^r_{1,1}\langle e_s, \widehat{\xi}^r_{1,1}\rangle\langle\widehat{\xi}^r_{1,1}, e_t\rangle)\mathbf{1}_{\Lambda_4}(i_1,r_1,j_1,l_1)\\
&\quad + (e_s{}^*e_t)\mathbf{1}_{\Lambda_1}(i_1,r_1,j_1,l_1) + \langle e_s, e_t\rangle\mathbf{I}_d\mathbf{1}_{\Lambda_2}(i_1,r_1,j_1,l_1)\\
&\quad + (e_t{}^*e_s)\mathbf{1}_{\Lambda_3}(i_1,r_1,j_1,l_1).
\end{aligned}$$

The discussion above allows us to affirm the following statements:

if $i = r = j = l$,

$$\langle\widetilde{H}^{i,r}e_s, \widetilde{H}^{j,l}e_t\rangle_n$$

(B.7a)
$$\begin{aligned}
&= \frac{1}{u(r)^2}\sum_{k=1}^n \rho^{-2k}\mathbf{X}_{k-1}(r)\mathbb{E}(\widehat{\xi}^r_{1,1} {}^*\widehat{\xi}^r_{1,1}{}^*e_s\widehat{\xi}^r_{1,1} {}^*\widehat{\xi}^r_{1,1}e_t)\\
&\quad + \frac{1}{u(r)^2}\sum_{k=1}^n \rho^{-2k}\mathbf{X}_{k-1}(r)(\mathbf{X}_{k-1}(r)-1)\{{}^*e_s e_t\mathbf{I}_d + e_t{}^*e_s + e_s{}^*e_t\}\\
&\quad - \frac{1}{u(r)^2}\sum_{k=1}^n \rho^{-2k}\mathbf{X}_{k-1}(r)^2\{e_s{}^*e_t\};
\end{aligned}$$

if $i = j \ne r = l$,

(B.7b)
$$\langle\widetilde{H}^{i,r}e_s, \widetilde{H}^{j,l}e_t\rangle_n = \frac{1}{u(i)u(r)}\sum_{k=1}^n \rho^{-2k}\mathbf{X}_{k-1}(i)\mathbf{X}_{k-1}(r){}^*e_s e_t\mathbf{I}_d;$$

if $i = l \ne j = r$,

(B.7c)
$$\langle \widetilde{H}^{i,r} e_s, \widetilde{H}^{j,l} e_t \rangle_n = \frac{1}{u(l)u(r)} \sum_{k=1}^{n} \rho^{-2k} \mathbf{X}_{k-1}(l) \mathbf{X}_{k-1}(r) e_t{}^* e_s;$$

(B.7d)  if $(j,l) \notin \{(i,r),(r,i)\}$,        $\langle \widetilde{H}^{i,r} e_s, \widetilde{H}^{j,l} e_t \rangle_n = 0.$

From properties (B.4) and (B.7a–d) it follows that, a.s. on $\mathbf{E}$,

$$\frac{1}{n}\langle \widetilde{H}^{i,r} e_s, \widetilde{H}^{j,l} e_t \rangle_n \underset{n \to \infty}{\longrightarrow} \frac{\mathbf{W}^2}{\rho^2}({}^*e_s e_t \mathbf{I}_d + e_t{}^* e_s) \qquad \text{if } i = r = j = l,$$

$$\frac{1}{n}\langle \widetilde{H}^{i,r} e_s, \widetilde{H}^{j,l} e_t \rangle_n \underset{n \to \infty}{\longrightarrow} \frac{\mathbf{W}^2}{\rho^2}{}^*e_s e_t \mathbf{I}_d \qquad\qquad \text{if } i = j \neq r = l,$$

$$\frac{1}{n}\langle \widetilde{H}^{i,r} e_s, \widetilde{H}^{j,l} e_t \rangle_n \underset{n \to \infty}{\longrightarrow} \frac{\mathbf{W}^2}{\rho^2} e_t{}^* e_s \qquad\qquad \text{if } i = l \neq r = j,$$

and

$$\frac{1}{n}\langle \widetilde{H}^{i,r} e_s, \widetilde{H}^{j,l} e_t \rangle_n \underset{n \to \infty}{\longrightarrow} 0 \qquad \text{if } (j,l) \notin \{(i,r),(r,i)\}.$$

Consequently,

$$\frac{1}{n}\langle \widetilde{H}^{i,r} e_s, \widetilde{H}^{j,l} e_t \rangle_n$$
$$\underset{n \to \infty}{\longrightarrow} \frac{\mathbf{W}^2}{\rho^2}\langle e_r, e_l \rangle \langle e_i, e_r \rangle \langle e_i, e_j \rangle({}^*e_s e_t \mathbf{I}_d + e_t{}^* e_s)$$
$$+ \frac{\mathbf{W}^2}{\rho^2}(\langle e_r, e_l \rangle \langle e_i, e_j \rangle(1 - \langle e_i, e_r \rangle)({}^*e_s e_t \mathbf{I}_d)$$
$$+ \langle e_i, e_l \rangle \langle e_r, e_j \rangle(1 - \langle e_i, e_r \rangle)e_t{}^* e_s)$$
$$= \frac{\mathbf{W}^2}{\rho^2}({}^*e_i \mathbf{I}_d e_j \langle e_r, e_l \rangle \langle e_s, e_t \rangle \mathbf{I}_d + {}^*e_i(e_l{}^* e_r)e_j(e_t{}^* e_s)) \qquad \text{a.s. on } \mathbf{E},$$

so

$$\frac{1}{n}\langle \widetilde{\mathsf{H}}^p, \widetilde{\mathsf{H}}^q \rangle_n \underset{n \to \infty}{\longrightarrow} \frac{\mathbf{W}^2}{\rho^2}(\langle e_r, e_l \rangle \langle e_s, e_t \rangle \mathbf{I}_d \otimes \mathbf{I}_d + (e_l{}^* e_r) \otimes (e_t{}^* e_s)) \qquad \text{a.s. on } \mathbf{E}.$$

Finally,

$$\frac{1}{n}\langle \widetilde{\mathcal{H}} \rangle_n \underset{n \to \infty}{\longrightarrow} \frac{\mathbf{W}^2}{\rho^2}((\mathbf{I}_d \otimes \mathbf{I}_d) \otimes (\mathbf{I}_d \otimes \mathbf{I}_d) + \mathcal{J} \boxtimes \mathcal{J}) \qquad \text{a.s. on } \mathbf{E},$$

where $\mathcal{J}$ is defined by (B.5).

(C) *Behavior of* $(\langle \mathcal{H}, \widetilde{\mathcal{H}} \rangle_n)$. For $1 \leq i, j \leq d, 1 \leq r, l \leq d$ and $1 \leq s, t \leq d,$

$\langle H^{i,r} e_s, \widetilde{H}^{j,l} e_t \rangle_n$

$$= \sum_{k=1}^{n} \rho^{-2k} \mathbf{M}_{k-1}^i \mathbb{E}(^*(\Delta \mathbf{M}_k^r) \langle \Delta \mathbf{M}_k^j, e_s \rangle \langle \Delta \mathbf{M}_k^l, e_t \rangle | \mathcal{G}_{k-1})$$

$$+ \sum_{k=1}^{n} \rho^{-2k} \langle \mathbf{M}_{k-1}^j, e_s \rangle \mathbb{E}(\Delta \mathbf{M}_k^{i*}(\Delta \mathbf{M}_k^r) \langle \Delta \mathbf{M}_k^l, e_t \rangle | \mathcal{G}_{k-1})$$

$$= \frac{1}{\sqrt{u(r)u(j)u(l)}}$$

$$\times \sum_{k=1}^{n} \rho^{-2k} \mathbf{M}_{k-1}^i \mathbb{E}(^*\zeta_k^r (K^r)^{-1/2*} e_s (K^j)^{-1/2} \zeta_k^{j*} \zeta_k^l (K^l)^{-1/2} e_t | \mathcal{G}_{k-1})$$

$$+ \frac{1}{\sqrt{u(r)u(i)u(l)}}$$

$$\times \sum_{k=1}^{n} \rho^{-2k} \langle \mathbf{M}_{k-1}^j, e_s \rangle \mathbb{E}((K^i)^{-1/2} \zeta_k^{i*} \zeta_k^r (K^r)^{-1/2*} \zeta_k^l (K^l)^{-1/2} e_t | \mathcal{G}_{k-1});$$

hence,

$$\langle H^{i,r} e_s, \widetilde{H}^{j,l} e_t \rangle_n$$

(B.8)
$$= \frac{1}{u(r)^{3/2}} \sum_{k=1}^{n} \rho^{-2k} \mathbf{X}_{k-1}(r) \mathbf{M}_{k-1}^i \mathbb{E}(^*\widehat{\xi}_{1,1}^r {}^* e_s \widehat{\xi}_{1,1}^r {}^* \widehat{\xi}_{1,1}^r e_t) \mathbf{1}_{\{j=l=r\}}$$

$$\times \frac{1}{u(r)^{3/2}} \sum_{k=1}^{n} \rho^{-2k} \mathbf{X}_{k-1}(r)^* e_s \mathbf{M}_{k-1}^j$$

$$\times \mathbb{E}((K^i)^{-1/2} \widehat{\xi}_{1,1}^i {}^* \widehat{\xi}_{1,1}^r (K^r)^{-1/2*} \widehat{\xi}_{1,1}^l (K^l)^{-1/2} e_t) \mathbf{1}_{\{i=l=r\}}.$$

Consequently,

$$\frac{1}{n} \langle H^{i,r} e_s, \widetilde{H}^{j,l} e_t \rangle_n \underset{n \to \infty}{\longrightarrow} 0 \qquad \text{a.s. on } \mathbf{E} \quad \text{and} \quad \frac{1}{n} \langle \mathcal{H}, \widetilde{\mathcal{H}} \rangle_n \underset{n \to \infty}{\longrightarrow} 0 \qquad \text{a.s. on } \mathbf{E}.$$

STEP 3—Verification of the Lindeberg condition for the martingale $\mathbb{H}$. For all $\varepsilon > 0$, we have

$$\frac{1}{n} \sum_{k=1}^{n} \mathbb{E}(\|\Delta \mathbb{H}_k\|^2 \mathbf{1}_{\{\|\mathbb{H}_k\| > \varepsilon \sqrt{n}\}} | \mathcal{G}_{k-1}) \underset{n \to \infty}{\longrightarrow} 0 \qquad \text{a.s. on } \mathbf{E}.$$

In fact, thanks to the properties

$$\max_{1 \le k \le n} \|\Delta \mathcal{H}_k\| = O(\ln n) \quad \text{and} \quad \max_{1 \le k \le n} \|\Delta \widetilde{\mathcal{H}}_k\| = O(\ln n),$$

which hold a.s. on $\mathbf{E}$, we have

$$\mathbf{1}_{\{\|\mathcal{H}_{k+1}\| + \|\widetilde{\mathcal{H}}_{k+1}\| > \varepsilon \sqrt{n}\}} \underset{n \to \infty}{\longrightarrow} 0 \qquad \text{a.s. on } \mathbf{E}.$$

Consequently,

$$\frac{1}{n}\sum_{k=1}^{n}\mathbb{E}(\|\Delta\mathbb{H}_k\|^2\mathbf{1}_{\{\|\mathbb{H}_k\|>\varepsilon\sqrt{n}\}}|\mathcal{G}_{k-1})$$

$$\leq\frac{1}{n}\sum_{k=1}^{n}\mathbb{E}(\|\Delta\mathcal{H}_k\|^2\mathbf{1}_{\{\|\mathcal{H}_k\|+\|\widetilde{\mathcal{H}}_k\|>\varepsilon\sqrt{n}\}}|\mathcal{G}_{k-1})$$

$$+\frac{1}{n}\sum_{k=1}^{n}\mathbb{E}(\|\Delta\widetilde{\mathcal{H}}_k\|^2\mathbf{1}_{\{\|\mathcal{H}_k\|+\|\widetilde{\mathcal{H}}_k\|>\varepsilon\sqrt{n}\}}|\mathcal{G}_{k-1})$$

$$=o\left(\frac{1}{n}\sum_{k=1}^{n}\rho^{-k}\|\mathcal{M}_k\|^2\right)+o\left(\frac{1}{n}\sum_{k=1}^{n}\mathbb{E}(\rho^{-2k}\|\Delta\mathcal{M}_k\|^4|\mathcal{G}_{k-1})\right)$$

$$=o(1)\qquad\text{a.s. on }\mathbf{E},$$

because the r.v.'s

$$\left(\frac{1}{n}\sum_{k=1}^{n}\rho^{-k}\|\mathcal{M}_k\|^2\right)\quad\text{and}\quad\left(\frac{1}{n}\sum_{k=1}^{n}\mathbb{E}(\rho^{-2k}\|\Delta\mathcal{M}_k\|^4|\mathcal{G}_{k-1})\right)$$

are almost-surely bounded on $\mathbf{E}$. The Lindeberg condition for $(\mathbb{H}_n)$ is proved.

STEP 4—Weak limit of $(\frac{1}{\sqrt{n}}\mathbb{H}_n)$.   From the classical martingale central limit theorem [17, 27], it follows that, conditional on $\mathbf{E}$,

$$(\text{B.9})\qquad\frac{1}{\sqrt{n}}\mathbb{H}_n\xrightarrow[n\to\infty]{\mathcal{L}}\begin{pmatrix}\Sigma_1\left(\left(\frac{2\mathbf{W}^2}{\rho^2(\rho-1)}\right)^{1/4}\mathbf{I}_d\right)\\[2mm]\Sigma_2\left(\sqrt{\dfrac{\mathbf{W}}{\rho}}\mathbf{I}_d\right)\end{pmatrix},$$

where $(\Sigma_r(T))_{r\in\{1,2\}}$ are independent identically distributed Gaussian vectors, which are also independent of the r.v. $\mathbf{W}$. The covariance matrix of $\Sigma_r(T)$ is equal to $(T\otimes T)\otimes(T\otimes Y)+{}^\perp(\text{Vect }T^*\text{Vect }T)\boxtimes{}^\perp(\text{Vect }T^*\text{Vect }T)$. Conditional on $\mathbf{E}_n$, this result remains true.   $\square$

PROOFS OF $(\mathrm{P}_9)$ AND $(\mathrm{P}_9')$.   Hereafter, we use Theorem 3 of [27] to prove a CLT for the couple of martingales $(\mathbb{H}_n,\mathcal{M}_n)$. Properties $(\mathrm{P}_9)$ and $(\mathrm{P}_9')$ are consequences of this result.

For the vectors $\widetilde{x}=(x_{(r-1)d+s})_{1\leq r,s\leq d}$, $\widetilde{z}=(z_{(l-1)d+t})_{1\leq l,t\leq d}$ and $y$ of $\mathbb{R}^{d^2}$, we set

$$\mathbb{H}_n'=\sum_{r=1}^{d}\sum_{s=1}^{d}\langle x_{(r-1)d+s},\mathsf{H}_n^{(r-1)d+s}\rangle+\sum_{l=1}^{d}\sum_{t=1}^{d}\langle z_{(l-1)d+t},\widetilde{\mathsf{H}}_n^{(l-1)d+t}\rangle,$$

$$\mathbb{U}_k^{(n)}=\frac{1}{\sqrt{n}}\Delta\mathbb{H}_k'+\langle y,\rho^{-n/2}\Delta\mathcal{M}_k\rangle,$$

$$\Xi_n(\widetilde{x}, \widetilde{z}, y) = \Xi_n((x_{(r-1)d+s})_{1 \le r, s \le d}, (z_{(l-1)d+t})_{1 \le l, t \le d}, y)$$

$$= \prod_{k=1}^{n} \mathbb{E}\{\exp\{i\Delta \mathbb{U}_k^{(n)}\} | \mathcal{G}_{k-1}\},$$

$$\Phi_n(y) = \prod_{k=1}^{n} \mathbb{E}(\exp\{i\langle y, \rho^{n/2} \Delta \mathcal{M}_k\rangle\} | \mathcal{G}_k),$$

$$\Psi_n(\widetilde{x}, \widetilde{z}) = \prod_{k=1}^{n} \mathbb{E}\left(\exp\left(i\frac{1}{\sqrt{n}}\Delta \mathbb{H}'_k\right) \Big| \mathcal{G}_{k-1}\right).$$

From the results proved in Steps 1, 2, 3 and 4, Case 3, the following classical property holds:

(B.10) $$\Psi_n(\widetilde{x}, \widetilde{z}) \xrightarrow[n \to \infty]{} \Psi_\infty(\widetilde{x}, \widetilde{z}) \qquad \text{a.s. on } \mathbf{E},$$

where

$$\Psi_\infty(\widetilde{x}, \widetilde{z})$$

$$= \exp\Bigg\{ -\frac{1}{2}\frac{2\mathbf{W}^2}{\rho^2(\rho-1)}$$

$$\times \Bigg( \sum_{r=1}^{d}\sum_{s=1}^{d}\sum_{l=1}^{d}\sum_{t=1}^{d} {}^* x_{(r-1)p+s}(\langle e_r, e_l\rangle\langle e_s, e_t\rangle \mathbf{I}_d \otimes \mathbf{I}_d$$

$$+ (e_l{}^*e_r) \otimes (e_t{}^*e_s))x_{(l-1)p+t}\Bigg)\Bigg\}$$

(B.11)
$$\times \exp\Bigg\{ -\frac{1}{2}\frac{\mathbf{W}^2}{\rho^2}$$

$$\times \Bigg( \sum_{r=1}^{d}\sum_{s=1}^{d}\sum_{l=1}^{d}\sum_{t=1}^{d} {}^* z_{(r-1)p+s}$$

$$\times (\langle e_r, e_l\rangle\langle e_s, e_t\rangle \mathbf{I}_d \otimes \mathbf{I}_d$$

$$+ (e_l{}^*e_r) \otimes (e_t{}^*e_s))z_{(l-1)p+t}\Bigg)\Bigg\}.$$

We also recall property (P$_4$), which states that

(B.12) $$\Phi_n(y) \xrightarrow[n \to \infty]{} \Phi_\infty(y) = \exp\left\{ -\frac{1}{2}\frac{\mathbf{W}}{(\rho-1)}\|y\|^2 \right\} \qquad \text{a.s. on } \mathbf{E}.$$

Now let us prove that, almost surely on $\mathbf{E}$,

(B.13) $$\mathbf{R}_n = |\Xi_n(\widetilde{x}, \widetilde{z}, y) - \Phi_n(y)\Psi_n(\widetilde{x}, \widetilde{z})| \xrightarrow[n \to \infty]{} 0.$$

For this purpose, the following inequalities, also used in [28], are relevant:

$$\mathbf{R}_n \leq \frac{1}{\sqrt{n}} \sum_{k=1}^{n} \mathbb{E}(|\langle y, \rho^{-n/2}\Delta\mathcal{M}_k\rangle| \times |\Delta\mathbb{H}'_k||\mathcal{G}_{k-1})$$

$$+ \frac{1}{4n} \sum_{k=1}^{n} \mathbb{E}(\langle y, \rho^{-n/2}\Delta\mathcal{M}_k\rangle^2|\mathcal{G}_{k-1}) \mathbb{E}((\Delta\mathbb{H}'_k)^2|\mathcal{G}_{k-1})$$

$$\leq \frac{1}{\sqrt{n}} \sum_{k=1}^{n} \mathbb{E}(\langle y, \rho^{-n/2}\Delta\mathcal{M}_k\rangle^2|\mathcal{G}_{k-1})^{1/2} \mathbb{E}((\Delta\mathbb{H}'_k)^2|\mathcal{G}_{k-1})^{1/2}$$

(B.14)
$$+ \frac{1}{4n} \sum_{k=1}^{n} \mathbb{E}(\langle y, \rho^{-n/2}\Delta\mathcal{M}_k\rangle^2|\mathcal{G}_{k-1}) \mathbb{E}((\Delta\mathbb{H}'_k)^2|\mathcal{G}_{k-1})$$

$$\leq \left(\frac{1}{n} \max_{1\leq k\leq n} \mathbb{E}((\Delta\mathbb{H}'_k)^2|\mathcal{G}_{k-1})\right)^{1/2} \left(\sum_{k=1}^{n} \rho^{-n}{}^*y\Delta\langle\mathcal{M}\rangle_k y\right)^{1/2}$$

$$+ \frac{1}{4n} \max_{1\leq k\leq n} \mathbb{E}((\Delta\mathbb{H}'_k)^2|\mathcal{G}_{k-1}) \sum_{k=1}^{n} \rho^{-n}{}^*y\Delta\langle\mathcal{M}\rangle_k y.$$

In fact, thanks to the inequality

$$\left|\prod_{k=0}^{n} a_k - \prod_{k=0}^{n} b_k\right| \leq \sum_{k=0}^{n} |b_k - a_k|,$$

which holds for $|a_k| \leq 1$ and $|b_k| \leq 1$, we have

$$\mathbf{R}_n \leq \sum_{k=1}^{n} \left|\mathbb{E}\left((1-\exp\{i\langle y, \rho^{-n/2}\Delta\mathcal{M}_k\rangle\})\left(1-\exp\left\{i\frac{1}{\sqrt{n}}\Delta\mathbb{H}'_k\right\}\right)\Big|\mathcal{G}_{k-1}\right)\right.$$

$$- \mathbb{E}((1-\exp\{i\langle y, \rho^{-n/2}\Delta\mathcal{M}_k\rangle\})|\mathcal{G}_{k-1})$$

$$\left. \times \mathbb{E}\left(\left(1-\exp\left\{i\frac{1}{\sqrt{n}}\Delta\mathbb{H}'_k\right\}\right)\Big|\mathcal{G}_{k-1}\right)\right|,$$

and the last quantity is less than the right-hand side of (B.14), thanks to the inequalities

$$\forall x \in \mathbb{R} \qquad |e^{ix}-1| \leq |x| \quad \text{and} \quad |e^{ix}-1-ix| \leq \frac{x^2}{2},$$

combined with the fact that $\mathbb{E}(\Delta\mathcal{M}_k|\mathcal{F}_{k-1}) = \mathbb{E}(\Delta\mathbb{H}'_k|\mathcal{F}_{k-1}) = 0$.

Because the sequence $(n^{-1}\langle\mathbb{H}'\rangle_n)$ is almost-surely convergent on $\mathbf{E}$, then

$$\frac{1}{n} \max_{1\leq k\leq n} \mathbb{E}((\Delta\mathbb{H}'_k)^2|\mathcal{G}_{k-1}) \underset{n\to\infty}{\longrightarrow} 0 \qquad \text{a.s. on } \mathbf{E}.$$

We also have

$$\sum_{k=1}^{n} (\rho^{-n}\,{}^*y\Delta\langle\mathcal{M}\rangle_k y) = O\left(\sum_{k=1}^{n}\rho^{-(n-k)/2}\right) = O(1) \qquad \text{a.s. on } \mathbf{E}.$$

Hence, property (B.13) is proved. Therefore, from (B.10) and (B.12),

$$\Xi_n(\tilde{x},\tilde{z},y)\underset{n\to\infty}{\longrightarrow}\Xi_\infty(\tilde{x},\tilde{z},y) = \Psi_\infty(\tilde{x},\tilde{z})\Phi_\infty(y) \qquad \text{a.s. on } \mathbf{E}.$$

Consequently, by Theorem 3 of [27], we can affirm that, conditional on $\mathbf{E}$,

$$\left\{\frac{1}{\sqrt{n}}\mathbb{H}_n, \rho^{-n/2}\mathcal{M}_n\right\}$$

(B.15)

$$\underset{n\to\infty}{\overset{\mathcal{L}}{\longrightarrow}}\left\{\begin{pmatrix}\Sigma_1\left(\left(\frac{2\mathbf{W}^2}{\rho^2(\rho-1)}\right)^{1/4}\mathbf{I}_d\right)\\\Sigma_2\left(\sqrt{\frac{\mathbf{W}}{\rho}}\mathbf{I}_d\right)\end{pmatrix}, \Sigma\left(\left(\frac{\mathbf{W}}{\rho-1}\right)^{1/2}\mathcal{I}_{d^2}\right)\right\},$$

where $\Sigma(\cdot)$ is a r.v. as in property (P5), and $(\Sigma_r(\cdot))_{r\in\{1,2\}}$ are r.v.'s independent of the pair $(\mathbf{W},\Sigma(\cdot))$ and distributed as in (B.9). Since property (B.15) is also valid conditional on $\mathbf{E}_n$, property (P9) for the martingale $(\mathcal{M}_n)$ follows from (B.15).

According to (B.2), (B.15) implies that, conditional on $\mathbf{E}$ or $\mathbf{E}_n$,

$$\left\{\left[\frac{1}{\sqrt{n}}\sum_{k=1}^{n}\rho^{-k}(\mathbf{M}_k^i{}^*\mathbf{M}_k^i - \langle\mathbf{M}^i\rangle_k)\right]_{1\le i\le d}, \rho^{-n/2}\mathcal{M}_n\right\}$$

$$\underset{n\to\infty}{\overset{\mathcal{L}}{\longrightarrow}}\left\{\left[\Sigma_1^i\left(\frac{\sqrt{2}\mathbf{W}}{(\rho-1)^{3/2}}\mathbf{I}_d\right) + \Sigma_2^i\left(\frac{\mathbf{W}}{\rho-1}\mathbf{I}_d\right)\right]_{1\le i\le d}, \Sigma\left(\left(\frac{\mathbf{W}}{\rho-1}\right)^{1/2}\mathcal{I}_{d^2}\right)\right\},$$

where $\{\Sigma_r^i(T); 1\le i\le d, r\in\{1,2\}\}$ are independent identically distributed Gaussian matrices, which we can choose independently of the pair $(\mathbf{W},\Sigma(\cdot))$. Moreover, the covariance matrix of $\mathrm{Vect}(\Sigma_r^i(T))$ is equal to $T\otimes T + {}^\perp(\mathrm{Vect}\,T\times{}^*\mathrm{Vect}\,T)$. Hence, (P$_9'$) is proved for the martingale $(\mathcal{M}_n)$. $\square$

PROOF OF PROPERTY (P10). Property (B.15) combined with

$$\sum_{k=1}^{n}\rho^{-k}\,\mathrm{tr}(\mathcal{M}_k{}^*\mathcal{M}_k - \langle\mathcal{M}\rangle_k)$$

(B.16)

$$= \left(\frac{\rho}{\rho-1}\right)\left\langle\begin{pmatrix}\mathrm{Vect}\,\mathcal{I}_{d^2}\\\mathrm{Vect}\,\mathcal{I}_{d^2}\end{pmatrix}, \mathbb{H}_n\right\rangle + o(\sqrt{n})$$

$$= \left(\frac{\rho}{\rho-1}\right)\left(\sum_{r=1}^{d}\sum_{s=1}^{d}{}^*e_s H_n^{r,r}e_s + \sum_{r=1}^{d}\sum_{s=1}^{d}{}^*e_s\widetilde{H}_n^{r,r}e_s\right) \qquad \text{a.s. on } \mathbf{E},$$

allows us to affirm that, conditional on $\mathbf{E}$ or $\mathbf{E}_n$,

$$\left\{ \frac{1}{\sqrt{n}} \sum_{k=1}^{n} \rho^{-k} \operatorname{tr}(\mathcal{M}_k{}^* \mathcal{M}_k - \langle \mathcal{M} \rangle_k), \rho^{-n/2} \mathcal{M}_n \right\} \xrightarrow[n \to \infty]{\mathcal{L}} \left\{ \tau G, \Sigma \left( \left( \frac{\mathbf{W}}{\rho - 1} \right)^{1/2} \mathcal{I}_{d^2} \right) \right\},$$

where $G$ is a standard Gaussian r.v., which is also independent of the pair $(\mathbf{W}, \Sigma(\cdot))$ and $\tau^2 = ((2\mathbf{W}^2 d^2 (\rho + 1))/(\rho - 1)^3$. In fact, according to (B.10) and (B.16), if $\tilde{x}$ and $\tilde{z}$ are the vectors of the canonical basis of $\mathbb{R}^{d^2}$, then

$$\begin{aligned}
\left( \frac{\rho - 1}{\rho} \right)^2 \tau^2 &= \frac{2\mathbf{W}^2}{\rho^2 (\rho - 1)} \sum_{r=1}^{d} \sum_{s=1}^{d} (\langle e_r, e_r \rangle^2 \langle e_s, e_s \rangle^2 + \langle e_r, e_r \rangle^2 \langle e_s, e_s \rangle^2) \\
&\quad + \frac{\mathbf{W}^2}{\rho^2} \sum_{r=1}^{d} \sum_{s=1}^{d} (\langle e_r, e_r \rangle^2 \langle e_s, e_s \rangle^2 + \langle e_r, e_r \rangle^2 \langle e_s, e_s \rangle^2) \\
&= \frac{4\mathbf{W}^2 d^2}{\rho^2 (\rho - 1)} + \frac{2\mathbf{W}^2 d^2}{\rho^2} = \frac{2\mathbf{W}^2 d^2 (\rho + 1)}{\rho^2 (\rho - 1)}.
\end{aligned}$$

The proof of Lemma 4.4 is complete. □

PROOF OF PROPERTY $(P_8)$. Let us prove the LIL announced in Lemma 3.4 for $(\mathcal{M}_n)$. The martingale

$$\overline{\mathcal{L}}_n = \left\langle \begin{pmatrix} \operatorname{Vect}(\mathcal{I}_{d^2}) \\ \operatorname{Vect}(\mathcal{I}_{d^2}) \end{pmatrix}, \mathbb{H}_n \right\rangle$$

satisfies the LIL

$$\limsup_{n \to \infty} (2\langle \overline{\mathcal{L}} \rangle_n \ln \ln \langle \overline{\mathcal{L}} \rangle_n)^{-1/2} |\overline{\mathcal{L}}_n| = 1 \qquad \text{a.s. on } \mathbf{E}.$$

This fact is a consequence of the LIL stated in [6, 7]. In fact, the properties

$$\| \rho^{-n/2} \mathcal{M}_n \|^2 = O(\ln n), \qquad \langle \overline{\mathcal{L}} \rangle_n = O(n)$$

and

$$\mathbb{E}(\| \rho^{-n/2} \Delta \mathcal{M}_n \|^4 | \mathcal{G}_{n-1}) = O((\ln n)^2),$$

a.s. on $\mathbf{E}$, imply that $\sum_{n \geq 1} \mathbb{E}(|(\langle \overline{\mathcal{L}} \rangle_n)^{-1/2} \Delta \overline{\mathcal{L}}_n|^4 | \mathcal{G}_{n-1}) < \infty$ a.s. on $\mathbf{E}$. Since

$$n^{-1} \langle \overline{\mathcal{L}} \rangle_n \xrightarrow[n \to \infty]{} \frac{2\mathbf{W}^2 d^2 (\rho + 1)}{\rho^2 (\rho - 1)} \qquad \text{a.s. on } \mathbf{E},$$

we deduce the LIL

$$\limsup \left( \frac{n}{\ln \ln n} \right)^{1/2} \left| \frac{1}{n} \sum_{k=1}^{n} \rho^{-k} \operatorname{tr}(\mathcal{M}_k{}^* \mathcal{M}_k - \langle \mathcal{M} \rangle_k) \right| = 2d \frac{\mathbf{W}(\rho + 1)^{1/2}}{(\rho - 1)^{3/2}}$$

$$\text{a.s. on } \mathbf{E},$$

taking in account the relation (B.16).

The LIL announced in the second part of ($P_8$) may be deduced easily from the last property, since we have $\widetilde{\mathcal{H}}_n = \frac{\mathbf{W}}{\rho}\widetilde{\mathcal{Z}}_n + O(\ln n)$ a.s. on $\mathbf{E}$. $\quad\square$

**Acknowledgment.** The authors would like to thank the anonymous referee for his discerning suggestions and comments about the original version of this paper. A substantial number of the references were brought to our attention by him.

FACULTÉ DES SCIENCES DE TUNIS
DÉPARTEMENT DE MATHÉMATIQUES
CAMPUS UNIVERSITAIRE
2092–EL MANAR 1
TUNIS
TUNISIA
E-MAIL: Faiza.Maaouia@fst.rnu.tn

FACULTÉ DES SCIENCES DE BIZERTE
DÉPARTEMENT DE MATHÉMATIQUES
CAMPUS UNIVERSITAIRE
7021–JARZOUNA BIZERTE
TUNISIA
E-MAIL: Abder.Touati@fsb.rnu.tn